\documentclass[11pt,a4paper,reqno]{amsart}
\usepackage[marginratio=1:1,totalwidth=15.75cm,totalheight=22.275cm]{geometry}
\usepackage[utf8]{inputenc}
\usepackage{amsmath,amssymb,amsthm,enumerate,amscd}
\usepackage[shortlabels]{enumitem}
\usepackage{graphicx}
\usepackage{csquotes}

\usepackage[dvipsnames]{xcolor}
\usepackage[colorlinks=true,linkcolor=blue!60!black,citecolor=teal!80!black,urlcolor=RoyalBlue]{hyperref}
\usepackage{tikz-cd}
\usepackage{amscd}

\numberwithin{equation}{section}

\usepackage{soul}

\newcommand*{\defeq}{\mathrel{\vcenter{\baselineskip0.5ex \lineskiplimit0pt
                     \hbox{\scriptsize.}\hbox{\scriptsize.}}}%
                     =}

\newcommand{\N}{\mathbb{N}}

\newcommand{\C}{\mathbb{C}}

\newcommand{\D}{\mathbb{D}}
\newcommand{\la}{\lambda}

\def\Cx{\mathbb{C}}
\def\Chat{\widehat{\mathbb{C}}}

\renewcommand{\hat}{\widehat}

\DeclareMathOperator{\Mod}{Mod}
\DeclareMathOperator{\inj}{inj}

\newcommand{\DD}{\mathbb{D}}
\DeclareMathOperator{\Def}{Def}
\DeclareMathOperator{\QC}{QC}
\DeclareMathOperator{\moduli}{mod}
\DeclareMathOperator{\aut}{Aut}
\newcommand{\length}{\ell}

\newcommand{\Id}{\operatorname{Id}}

\DeclareMathOperator{\Per}{Per}
\DeclareMathOperator{\GO}{GO}

\newcommand{\T}{\mathcal{T}}

\def\C{{\mathbb{C}}}
\def\D{{\mathbb{D}}}
\def\N{{\mathbb{N}}}

\theoremstyle{plain}
\newtheorem{thm}{Theorem}[section]
\newtheorem{lemma}[thm]{Lemma}
\newtheorem{cor}[thm]{Corollary}

\newtheorem{obs}[thm]{Observation}

\newtheorem{question}{Question}

\theoremstyle{definition}
\newtheorem{defi}[thm]{Definition}

\newtheorem*{structure}{Structure of the article}

\theoremstyle{remark}

\newtheorem*{remark}{Remark}

\title[Teichm\"uller spaces and normal forms in wandering domains]{Teichm\"uller spaces and normal forms associated to wandering domains}

\author[N. Fagella]{Núria Fagella}
\author[G. R. Ferreira]{Gustavo R. Ferreira}
\author[L. Pardo-Simón]{Leticia Pardo-Simón}

\address{\noindent Dept. de Matemàtiques i Informàtica\\ Universitat de Barcelona\\ Catalonia\\ Spain\\
	\newline  Centre de Recerca Matemàtica\\ Bellaterra\\ Catalonia\\ Spain.
	\textsc{\newline \indent 
		\href{https://orcid.org/0000-0002-5466-0579%
		}{\includegraphics[width=1em,height=1em]{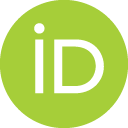} {\normalfont https://orcid.org/0000-0002-5466-0579}}
}}
\email{nfagella@ub.edu}

\address{\noindent Centre de Recerca Matemàtica\\ Bellaterra\\ Catalonia\\ Spain.\\
	\textsc{\newline \indent 
		\href{https://orcid.org/0000-0002-7330-0018%
		}{\includegraphics[width=1em,height=1em]{orcid2.png} {\normalfont https://orcid.org/0000-0002-7330-0018}}
}}
\email{grodrigues@crm.cat}
\address{\noindent Dept. de Matemàtiques i Informàtica\\ Universitat de Barcelona\\ Catalonia\\ Spain\\
	\newline  Centre de Recerca Matemàtica\\ Bellaterra\\ Catalonia\\ Spain.
	\textsc{\newline \indent 
		\href{https://orcid.org/0000-0003-4039-5556%
		}{\includegraphics[width=1em,height=1em]{orcid2.png} {\normalfont https://orcid.org/0000-0003-4039-5556}}
}}
\email{lpardosimon@ub.edu}

\thanks{This work was partially supported by the project PID2023-147252NB-I00 financed by MICIU/AEI MCIN/AEI/10.13039/501100011033,	FEDER, EU; by the Spanish State Research Agency, through the Severo Ochoa and María de Maeztu Program for Centers and Units of	Excellence in R\&D (CEX2020-001084-M); and by the Catalan government through ICREA Academia 2020. The second author acknowledges financial support from the European Union's Horizon Europe research and innovation programme under the Marie Skłodowska-Curie Grant Agreement No. 101208256. The third author is a Serra Húnter fellow.
}

\begin{document}
\begin{abstract}
We study the dynamical Teichm\"uller space ${\mathcal T}(U,f)$ associated to a wandering domain $U$ of an entire function $f$. We show that a discrete grand orbit relation in $U$ forces ${\mathcal T}(U,f)$ to be infinite dimensional, thereby answering a question of Fagella--Henriksen. We further describe the geometry of these spaces by developing normal forms for the dynamics on wandering domains, yielding global linearising coordinates in the discrete case and power-type dynamics between annuli in the indiscrete case.
\end{abstract}

\maketitle

\section{Introduction}

Let $f\colon S \to S$ be a holomorphic endomorphism of a Riemann surface $S$. The iterates
$f^n\colon S\to S$ of $f$ give rise to the classical partition of $S$, first observed by Fatou and Julia in the 1920s, into the \emph{Fatou set} $F(f)$ and the \emph{Julia set} $J(f)$. The Fatou set is the largest open set on which the family $(f^n)_n$ is normal; its complement $J(f)$ is the locus of chaotic dynamics. 

As in many areas of dynamics, one is often interested not only in a single map, but in how maps vary in parameter space. In complex dynamics this leads naturally to the study of moduli spaces of rational and entire maps. After a series of landmark developments \cite{MSS83,DH85}, McMullen and Sullivan \cite{MS98} introduced a highly geometric framework for this problem: the \emph{dynamical Teichmüller space} of a holomorphic map (see \S~\ref{ssec:teichmuller}). This point of view has had many applications, for instance to the study of conjugacy classes and structural stability of holomorphic dynamical systems, and it reveals a deep connection between complex dynamics and Teichmüller theory.

In order to describe our results, we begin by recalling some basic dynamical notions. Given a rational or entire map $f\colon S\to S$ (where $S$ is either the Riemann sphere or the complex plane, respectively), the \emph{grand orbit} of a point $z\in S$ is
\[
\GO(z) \defeq \{ w\in S : f^m(w) = f^n(z) \text{ for some } m,n\in\N\},
\]
that is, the set of points that can be connected to $z$ by forward or backward iteration of $f$. This defines an equivalence relation --- the \emph{grand orbit relation} --- whose equivalence classes are precisely the distinct grand orbits or, equivalently, the elements of $S/f$. We say that a set $V \subset S$ has a \emph{discrete grand orbit relation} if $\GO(z)$ is a discrete subset of $V$ for some (and hence for every, by \cite[Theorem~A]{EFGP24}) $z\in V$. Otherwise we say that the grand orbit relation in $V$ is \emph{indiscrete}. 

Following McMullen and Sullivan, we dynamically mark the set of singular values (critical values and asymptotic values), $S(f)$, and of periodic points, $\Per(f)$, of $f$, and consider the closure of their grand orbits:
\[
\hat J(f) \defeq \overline{\GO(S(f)\cup \Per(f))}, \qquad
\hat F(f) \defeq S\setminus \hat J(f),
\]
Then $\hat J(f)$ is completely invariant and contains $J(f)$, while $\hat F(f)$ is a completely invariant subset of the Fatou set which contains no points in the grand orbit of a singular value. In particular, the restriction $f\colon \hat F(f)\to \hat F(f)$ is a covering map. The connected components of $\hat F(f)$, or equivalently (abusing terminology) their grand orbits, will be called \emph{grand components}. 

The dynamical Teichmüller space $\T(f)$  supported on the Fatou set can be decomposed in terms of these grand components 
(see Theorem~\ref{thm:MSfactor} for a precise statement). For the purposes of this paper, it is convenient to write $\T(V,f)$ for the factor of $\T(f)$ associated to a given grand component $V$, which we call the {\em associated Teichm\"uller space} to $V$ (or to any component of~$V$).

When $f$ is a rational map, $\T(f)$ is mapped locally injectively into the (finite-dimensional) moduli space of $f$, and therefore $\T(f)$ itself is finite-dimensional (in fact, it was proved by Astorg \cite{Ast17} that this mapping is an immersion). In contrast, for entire functions $\T(f)$ may be infinite-dimensional: this happens, for instance, in attracting or parabolic basins with infinitely many singular values, in many Baker domains, and in most known examples of wandering domains; see \cite{FH06,FH09}. In all these examples the relevant factor of $\T(f)$ comes from a grand component with a discrete grand orbit relation, although there are also Baker domains with discrete grand orbit relation and finite-dimensional dynamical Teichmüller space. Additionally, Sullivan’s proof of the No Wandering Domains Theorem shows that a hypothetical wandering domain for a rational map would likewise have a \emph{discrete} grand orbit relation and would force an \emph{infinite}-dimensional dynamical Teichmüller space. This parallel between the rational case (hypothetically) and the known transcendental examples led Fagella and Henriksen to ask in \cite[§4.3]{FH09} the following.

\begin{question}[Fagella--Henriksen]\label{q:FH}
Can a wandering domain $U$ of an entire function $f$ have a discrete grand orbit relation and an associated finite-dimensional Teichmüller space?
\end{question}

Recent work shows that the situation is more subtle. In \cite[Theorem~B]{EFGP24} it is proved that discrete and indiscrete grand orbit relations may coexist inside a single wandering domain. More precisely, if $U$ is a wandering domain of an entire function $f$, the set $\hat F(f)\cap U$ may break into several connected components, some with discrete grand orbit relation and some with indiscrete grand orbit relation. This leads to the following refinement of Question~\ref{q:FH}.

\begin{question}\label{q:Q2}
	Let $U$ be a wandering domain of $f$. Can a connected component of $\hat F(f)\cap U$ have a discrete grand orbit relation and an associated finite-dimensional Teichmüller space?
\end{question}

In this paper we answer Question~\ref{q:Q2}, and hence Question~\ref{q:FH}, in the negative. More precisely, we prove the following.

\begin{thm}[Discreteness if and only if infinite-dimensional Teichmüller space]\label{thm:main}
	Let $f\colon \C\to\C$ be an entire function with a wandering domain $U$, let $V$ be a connected component of $\hat F(f)\cap U$, and let $\hat V$ denote the grand orbit of $V$. Then $\T(\hat V,f)$ is infinite-dimensional if and only if the grand orbit relation in $V$ is discrete.
\end{thm}
\begin{remark}
	By work of McMullen and Sullivan \cite{MS98} and Fagella and Henriksen \cite{FH09} (see Theorem \ref{thm:MSteich}), infinite-dimensional factors of $\T(f)$ associated to periodic components can only arise from grand components with discrete grand orbit relation. Theorem~\ref{thm:main} shows that, for wandering domains, this condition is also sufficient.
\end{remark}
Beyond this dichotomy, we develop a normal form theory for the internal dynamics on components of wandering domains. In the spirit of the classical normal forms for periodic Fatou components, we construct normalising coordinates on the connected components of $\hat F(f)\cap U$ and use them to give a more detailed description of the associated Teichmüller spaces, both in the discrete and in the indiscrete case.

To formulate our next result, it is convenient to use the concept of hyperbolic distortion. For a hyperbolic set $U\subset \C$ denote by $d_U$ the hyperbolic distance in $U$. If $f\colon U \to V$ is a holomorphic map between hyperbolic domains, the \emph{hyperbolic distortion} of $f$ at $z\in U$ is defined by 
$$ \|Df(z)\|_{U}^{V} := \lim_{z' \to z} \frac{d_V(f(z'), f(z))}{d_U(z',z)}, $$ and it equals the modulus of the hyperbolic derivative of $f$ at $z$.
 \begin{thm}[Discrete orbit relation and linearising coordinates] \label{thm:discrete_intro}
	Let $f\colon\Cx\to\Cx$ be an entire function with a wandering domain $U$, let $V$ be a component of $\widehat F(f)\cap U$ and let $V_n = f^n(V)$ for $n\geq 0$. Then, the following are equivalent.
 	\begin{enumerate}[label={\normalfont (\roman*)}]
 		\item \label{item:discrete}The grand orbit relation in $V$ is discrete; 
 		\item \label{item:linearising} For each $z\in V$, there exists a sequence of holomorphic maps $\phi_n\colon(V_n, f^n(z))\to(\Cx, 0)$, $n\geq 0$, and of neighbourhoods $W_n\subset V_n$ of $f^n(z)$  with $f(W_n)\subset W_{n+1}$, such that 
 		\[\phi_{n+1}\circ f\vert_{V_n} =  \lambda_n(z)\cdot \phi_n,
 		\]
 		where $\lambda_n(z)$ is the hyperbolic distortion $ \|Df(f^n(z))\|_{U_n}^{U_{n+1}}$, and so that $\phi_n$ maps $W_n$ biholomorphically onto a neighbourhood of zero.
 	\end{enumerate}
 \end{thm}

Theorem~\ref{thm:discrete_intro} is a global non-autonomous linearisation result for wandering domains, and may be viewed as an extension of results of Pommerenke \cite{Pom94} on non-autonomous dynamics on the disc (see §\ref{sec:proofteich}). We stress that any two sequences of nonzero linear maps $w\mapsto \lambda_n w$ and $w\mapsto \mu_n w$ are conformally conjugate by an affine change of coordinates on~$\C$. Thus the appearance of the hyperbolic distortions $\lambda_n(z)$ in Theorem~\ref{thm:discrete_intro} is not itself a new invariant; the substance of the theorem is that, inside a component of $\widehat F(f)\cap U$ with discrete grand orbit relation, the dynamics of $f$ can be globally conjugated to a sequence of linear maps.

As an application of Theorem~\ref{thm:discrete_intro}, we describe the
Teichmüller spaces associated with non-contracting simply connected wandering domains, i.e. whose
internal dynamics does not contract the hyperbolic metric too strongly along the
orbit (known as semi–contracting and
eventually isometric wandering domains); see \cite{BEFRS21} or Theorem~\ref{thm:Ainternal}.

\begin{thm}\label{thm:teich}
	Let $U$ be a simply connected wandering domain of the entire function $f$. Assume that for some (and hence for all) $z\in U$, there exists $c > 0$ such that $\|Df(f^{n}(z))\|_{U_{n}}^{U_{n+1}} \geq c$ for all $n$. Let $E = (\GO(S(f))\cap U)/f$ and $\hat U=U\cap \hat F(f)$. Then, the grand orbit relation in $U$ is discrete, and 
	\begin{enumerate}[label={\normalfont (\alph*)}]
		\item If $U$ is contracting, $\T(\hat U, f) \simeq \T(\Cx\setminus \tilde{E})$, and $\tilde{E}$ is at least countably infinite;\label{it:contr}
		\item If $U$ is not contracting, $\T(\hat U, f) \simeq \T(\DD\setminus \tilde{E})$,\label{it:notcontr}
	\end{enumerate}
	where $\tilde{E}$ is homeomorphic to $E$. In fact, there exists a biholomorphism $\phi\colon \GO(U)/f\to X$, where $X\in\{\DD, \Cx\}$, and $\tilde E = \phi(E)$.
\end{thm}

In the statement, $\T(\C\setminus \tilde{E})$ denotes the classical Teichm\"uller space of the Riemann surface $\C\setminus \tilde{E}$. Notice that the set $\tilde{E}$ may conceivably have non-empty interior, which would make the distinction between \ref{it:contr} and \ref{it:notcontr} irrelevant -- but in many cases, $\tilde{E}$ is in fact discrete.

The indiscrete counterpart of Theorem~\ref{thm:discrete_intro} also exhibits a strong rigidity phenomenon:

\begin{thm}[Indiscrete orbit relation and power maps]\label{thm:indiscrete_intro}
	Let $f\colon\Cx\to\Cx$ be an entire function with a wandering domain $U$, let
	$V$ be a component of $\widehat F(f)\cap U$, and denote $V_n = f^n(V)$ for
	$n\geq 0$. Then the following are equivalent:
	\begin{enumerate}[label={\normalfont (\roman*)}]
		\item the orbit relation in $V$ is indiscrete;
		\item $V$ is doubly connected and there is a sequence of conformal maps
		$\phi_n\colon V_n\to A_n$, $n\geq 0$, where $A_n$ is either an annulus of
		finite modulus or a punctured disc for all $n$, such that
		\[
		\phi_{n+1}(f(z)) = (\phi_n(z))^{d_n} \quad \text{for all } z\in V_n,
		\]
		where $d_n = \deg(f|_{V_n})$ and $d_n\geq 2$ for infinitely many values of $n$.
	\end{enumerate}
\end{thm}

This theorem can already be extracted from a careful reading of the proof of
\cite[Theorem~6.1]{MS98}; for completeness, we include here a more geometric
argument. The conclusion was also known for multiply connected wandering domains
of entire functions, which always have indiscrete grand orbit relation and
``power–map-like'' dynamics; see \cite[Corollary~6.9]{EFGP24} and \cite[Theorem 5.2]{BRS13}.

The combination of Theorem~\ref{thm:main} with Theorem~\ref{thm:indiscrete_intro}
gives a precise criterion for when the Teichmüller space associated with a
wandering domain is finite-dimensional.

\begin{cor}\label{cor:finite}
	Let $f\colon\Cx\to\Cx$ be an entire function with a wandering domain $U$. Then
	$\T(\widehat F(f)\cap U, f)$ is finite-dimensional if and only if there exists
	a finite number $N$ of grand orbits of components of $\widehat F(f)\cap U$,
	and each of these components has indiscrete orbit relation. In this case,
	$\widehat F(f)\cap U$ is a collection of $m$ annuli of finite modulus and
	$N-m$ punctured discs for some $0\leq m\leq N$, and $\dim\bigl(\T(\widehat F(f)\cap U, f)\bigr) = m.$
\end{cor}

An example of a wandering domain satisfying the properties in
Corollary~\ref{cor:finite} is given in \cite[p.~321, Example~2]{FH09} as the
logarithmic lift of a basin of attraction of a superattracting orbit in
$\Cx^*$. This might suggest that the existence of a wandering domain with
finite-dimensional Teichmüller space forces the presence of a
\emph{supercritical orbit}, i.e.\ an orbit containing infinitely many critical
points of $f$. However, this is not the case: by applying the procedure
described in \cite[Section~7]{EFGP24} to the Chebyshev polynomial
$p(z) = z^2 - 2$, one obtains an entire function $f$ with a wandering domain
$U$ whose grand orbit structure is reminiscent of that of a neighbourhood of
$J(p)$: $\hat J(f)\cap U$ is a continuum, and the grand orbit relation in
$\hat F(f)\cap U$ (which is an annulus of finite modulus) is indiscrete, but
there are no supercritical orbits at play.

\begin{remark}
	Our results concern only the factor of the dynamical Teichmüller space coming from the wandering domain inside $\widehat F(f)$, namely $\T(\widehat F(f)\cap U,f)$. In general, the complementary contribution supported on $\widehat J(f)\cap U$ may already be infinite-dimensional: the set $\widehat J(f)\cap U$ can have non-empty interior and support invariant line fields, leading to an infinite-dimensional space $\T(\widehat J(f)\cap U,f)$, and hence to an infinite-dimensional Teichmüller space $\T(f)$.
\end{remark}

\begin{structure} Section \ref{sec:prelims} has some preliminary results on Teichm\"uller theory and the dynamics of wandering domains (the expert reader is encouraged to skip this section and refer back to it as needed). Section \ref{sec:proofMT} contains a proof of Theorem \ref{thm:main}. In Sections \ref{sec:proofdiscrete}, we prove Theorem \ref{thm:discrete_intro}, followed by Theorem \ref{thm:teich} in Section \ref{sec:proofteich}. Finally, Theorem \ref{thm:indiscrete_intro} is proved in Section~\ref{sec:proofindiscrete}.
\end{structure}

\noindent\textbf{Acknowledgements.}
We thank Lukas Geyer for pointing out Sullivan's proof of Lemma \ref{lem:inflimits}, and Phil Rippon for reference \cite{Pom94}. We also thank the referee for their thorough comments.

\section{Preliminaries}\label{sec:prelims}

\subsection{Basic notation} We denote the closure of a set $A$ by $\overline{A}$. We write $A\subset B$ to indicate
that $A$ is a subset of $B$, including the possibility that $A=B$. 
For $r>0$, we denote by $\D_r(p)$  the open (Euclidean) disc
of radius $r$ centred at $p$. In the case $p=0$, we write $\D_r$ and if, additionally, $r=1$, we set $\D := \D_1$. For a
simply connected domain $S\subset \C$ or a hyperbolic Riemann surface $U$, $d_S(\cdot, \cdot )$ and  $\ell_S(\cdot)$ denote respectively
the hyperbolic distance and hyperbolic length in $S$, and $D_S(p, r)$ denotes the hyperbolic disc centred at $p$ of radius $r$. Sequences of points, sets, and functions such as
$(x_n)$ are indexed by $n\geq 0$, unless otherwise specified.

\subsection{The Teichm\"uller space of an entire function}\label{ssec:teichmuller}
The Teichm\"uller space of a holomorphic dynamical system was introduced by McMullen and Sullivan in \cite{MS98} (cf.\ \cite{EpsteinThesis}) as a tool to study the corresponding \emph{moduli space}. For our purposes, the moduli space of an entire function $f$ is
\[
\moduli(f) := \{g\in\mathcal{H}(\Cx)\colon \varphi\circ f = g\circ\varphi,\ \varphi\colon\Cx\to\Cx\ \text{quasiconformal}\}/\sim,
\]
where $g_1\sim g_2$ if $g_1$ and $g_2$ are conformally conjugate. We only recall the basic constructions needed later; for a detailed account in the entire setting see \cite{FH09,EpsteinThesis}.

Let $f\colon V\to V$ be a holomorphic endomorphism of a one-dimensional complex manifold\footnote{As usual, we reserve the term Riemann surface for \emph{connected} one-dimensional complex manifolds.}. If $g\colon V'\to V'$ is holomorphic and $\phi\colon V\to V'$ is quasiconformal with $g\circ\phi=\phi\circ f$, we write $(g,\phi)$ (note that $V'=\phi(V)$ is determined by $\phi$). Two such pairs $(g_1,\phi)$ and $(g_2,\psi)$ are equivalent if there exists a conformal map $c\colon\phi(V)\to\psi(V)$ such that $\psi = c\circ\phi$; we then write $(g_1,\phi)\sim (g_2,\psi)$. The \emph{deformation space} of $f$ is
\[
\Def(V,f) := \{(g,\phi)\colon \phi\colon V\to V' \text{ quasiconformal},\ g\colon V'\to V' \text{ holomorphic},\ g\circ\phi=\phi\circ f\}/\sim.
\]

If $\mu_\phi$ denotes the Beltrami form\footnote{As customary, we reserve the term \emph{Beltrami form} for measurable $(-1,1)$-tensors on $V$ with supremum norm strictly less than one. If $V\subset\Cx$, we also say \emph{Beltrami coefficient}.} of $\phi$, then by the Ahlfors--Bers theorem the map
$(g,\phi)\mapsto \mu_\phi$ induces a bijection between $\Def(V,f)$ and $M_1(V,f)$, the Banach manifold of $f$-invariant Beltrami forms on $V$ (see e.g.\ \cite[Theorem~4.4]{MS98}). In particular, $\Def(V,f)$ carries a natural Banach manifold structure.

Next, let $\QC(V,f)$ be the group of quasiconformal automorphisms of $V$ commuting with $f$, and let $\QC_0(V,f)$ be the normal subgroup of those isotopic to the identity \emph{relative the ideal boundary of $V$}, via an isotopy through elements of $\QC(V,f)$; see \cite{EM88} for details. The group $\QC(V,f)$ acts on $\Def(V,f)$ by
\[
(g,\phi)\mapsto (g,\phi\circ\omega^{-1}),\qquad \omega\in\QC(V,f).
\]

\begin{defi}
	The \emph{Teichm\"uller space} of the holomorphic map $f\colon V\to V$ is
	\[
	\T(V,f) := \Def(V,f)/\QC_0(V,f).
	\]
\end{defi}

An element of $\T(V,f)$ is thus represented by an equivalence class $[(g,\phi)]$, where we quotient first by postcomposition with conformal maps (the definition of $\Def(V,f)$) and then by $f$-equivariant isotopies relative the ideal boundary. Forgetting the marking $\phi$ gives a natural map $\T(V,f)\longrightarrow \moduli(f).$ When $f$ is rational, McMullen and Sullivan showed that the induced action of the mapping class group $\Mod(\Chat,f):=\QC(\Chat,f)/\QC_0(\Chat,f)$ on $\T(\Chat,f)$ is properly discontinuous \cite[Theorem~2.3]{MS98}; in particular, the quotient $\T(\Chat, f)/\Mod(\Chat, f) \simeq \moduli(f)$ inherits an orbifold structure in much the same way the quotient of a Riemann surface by a group of automorphisms acting freely and properly discontinuously inherits a complex structure. Under this orbifold structure, the Teichm\"uller space is realised as a universal orbifold cover of $\moduli(f)$. Notice, however, that this orbifold structure is not always compatible with the ambient structure of $\mathrm{rat}_d$ \cite[Section 10.2]{Bra93}.

We now specialise to entire functions. Let $f$ be transcendental entire. Recall from the introduction the dynamically marked sets $\hat J(f)$ and $\hat F(f) := \Cx\setminus\hat J(f)$; the map $f|_{\hat F(f)}$ is a holomorphic self-covering of $\hat F(f)$, and it decomposes into restrictions of the form
\[
f\colon \hat V \to \hat V,\qquad \hat V := \GO(\hat U),
\]
where $\hat U$ runs over the connected components of $\hat F(f)$. A fundamental structural result is the following factorisation of the dynamical Teichm\"uller space.

\begin{thm}[{\cite[Theorem~3.3]{FH09}, \cite{MS98}}]\label{thm:MSfactor}
	Let $f$ be a transcendental entire function. Let $\hat V_i$ denote the collection of pairwise disjoint grand orbits of connected components of $\hat F(f)$. Then
	\[
	\T(\Cx,f) \;=\; M_1(\hat J(f),f)\times \prod_i^* \T(\hat V_i,f),
	\]
	where $\prod^*$ denotes the restricted product of Teichm\"uller spaces (see \cite[p.~304]{FH09}).
\end{thm}

Thus, up to the space $M_1(\hat J(f),f)$ of invariant line fields supported on $\hat J(f)$, understanding $\T(\Cx,f)$ reduces to understanding the Teichm\"uller spaces $\T(\hat V_i,f)$ for grand orbits inside $\hat F(f)$. The geometry of these pieces is described by the following theorem.

\begin{thm}[{\cite[Theorem~6.1]{MS98}, \cite[Theorem~3.4]{FH09}}]\label{thm:MSteich}
	Let $f$ be a transcendental entire function, and let $\hat V\subset\hat F(f)$ be the grand orbit of a connected component of $\hat F(f)$. Then:
	\begin{enumerate}[label={\normalfont (\arabic*)}]
		\item If the grand orbit relation of $f|_{\hat V}$ is discrete, then $\hat V/f$ is a Riemann surface and
		\[
		\T(\hat V,f) \simeq \T(\hat V/f).
		\]
		\item If the grand orbit relation of $f|_{\hat V}$ is indiscrete and some (hence every) connected component of $\hat V$ is an annulus of finite modulus, then
		\[
		\T(\hat V,f) \simeq \mathbb{H}.
		\]
		\item Otherwise, $\T(\hat V,f)$ is trivial.
	\end{enumerate}
\end{thm}

Very briefly, in case (1) every $f$-invariant Beltrami coefficient on $\hat V$ descends to a Beltrami form on the quotient $\hat V/f$, so deformations of $(\hat V,f)$ are in bijection with deformations of the Riemann surface $\hat V/f$. In case (2), the grand orbit structure induces an $S^1$-foliation on an annular component, and one can deform either its modulus or the relative twisting of its boundary components, giving a one-dimensional Teichm\"uller space; see \cite[Subsection~5.1]{MS98}. In case (3) each component of $\hat V$ is a punctured disc, foliated again by $S^1$-orbits, and the only deformations are conformal.

\subsection{Wandering domains: internal dynamics and orbit relations}\label{subsec:wd}
To understand the Teichm\"uller space of a wandering domain, one must understand how the orbits of its points relate to one another. We start recalling the following classification result.
\begin{thm}[Classification of simply connected wandering domains, {\cite[Theorem~A]{BEFRS21}}]\label{thm:Ainternal}Let $U$ be a simply connected wandering domain of a transcendental entire function $f$ and let $U_n$ be the Fatou component containing $f^n(U)$, for $n \geq 0$. Define the countable set of pairs  
	$$ E = \{(z, z') \in U \times U : z\neq z' \text{ and } f^k(z) = f^k(z') \text{ for some } k \in \mathbb{N} \}. $$  
	Then, exactly one of the following holds.  
	\begin{enumerate}
		\item $ d_{U_n}(f^n(z), f^n(z'))\underset{n\to\infty}{\longrightarrow} c(z,z') = 0 $ for all $z, z' \in U$, and we say that $U$ is \emph{contracting}.  
		\item $ d_{U_n}(f^n(z), f^n(z')) \underset{n\to\infty}{\longrightarrow} c(z,z') > 0 $ and $ d_{U_n}(f^n(z), f^n(z')) \neq c(z, z') $ for all $(z, z') \in (U \times U) \setminus E$, $n \in \mathbb{N}$, and we say that $U$ is \emph{semi-contracting}.  
		\item There exists $N > 0$ such that for all $n \geq N$, $ d_{U_n}(f^n(z), f^n(z')) = c(z,z') > 0 $ for all $(z, z') \in (U \times U) \setminus E$, and we say that $U$ is \emph{eventually isometric}.  
	\end{enumerate}
	Moreover, for $z \in U$ and $n \in \mathbb{N}$ let $\lambda_n(z)$ be the hyperbolic distortion $\|Df(f^n(z))\|_{U_n}^{U_{n+1}}$. Then  
	\begin{itemize}
		\item $U$ is contracting if and only if $\sum_{n=0}^{\infty} (1 - \lambda_n(z)) = \infty$.  
		\item $U$ is eventually isometric if and only if $\lambda_n(z) = 1$ for $n$ sufficiently large.  
	\end{itemize}
\end{thm}

Next, let $f$ be an entire function with a wandering domain $U$,  and for each $n\in \N$, let $U_n$ be the Fatou component containing $f^n(U)$. By using Riemann maps, one can relate the dynamics of $f$ restricted to $(U_n)$ to a sequence of self-maps of the unit disc. These maps happen to be \textit{inner functions}, that is,  holomorphic self-maps of $\D$ preserving the unit circle Lebesgue-almost everywhere in the sense of radial limits (see e.g., \cite{associates_20}). More precisely, we fix a point $z_0\in U=U_0$ and for each $n\geq 0$ we choose a Riemann map $\psi_n\colon U_n\to \D$ such that $\psi_n(f^n(z_0))=0$. Then, we consider the sequence of holomorphic maps $g_n\colon \D\to\D$ given by 
\begin{equation}\label{eq:associated_inner}
	g_n := \psi_{n+1} \circ f \circ \psi_n^{-1},
\end{equation}
and call this an \textit{associated sequence} of inner functions.  Sometimes, it will be convenient to think of the domains and ranges of the maps $g_n$ as distinct copies of the unit disc, denoted $\D^{(n)}$, where
$g_n:\D^{(n)} \to \D^{(n+1)}$ for any $n\geq 0$.  It then holds, see \cite[Lemma 3.1]{BEFRS21}, that for each $n\in \N$,
\begin{equation}\label{eq_deriv0}
	\lambda_n(z_0):= \|Df(f^n(z))\|_{U_n}^{U_{n+1}}=\vert g'_n(0)\vert \leq 1,
\end{equation}
where the last inequality follows from Schwarz's lemma.

We will additionally use the fact that for $G_n := g_{n-1}\circ \cdots \circ g_0$, 
\begin{equation}\label{eq:Gnto0}
	\sum_{n=0}^{\infty} (1 - \lambda_n(z)) = \infty  \implies G_n(z)\to 0 \text{ as } n\to \infty, \text{ for all } z\in\D,
\end{equation}
as shown in \cite[Theorem 2.1]{BEFRS21}.

Finally, we include a result from \cite{EFGP24} on grand orbit relations on Fatou components, stated for the case that interests us, namely wandering orbits.
\begin{thm}[{\cite[Theorem~A]{EFGP24}}]\label{thm:GOrelations}
	Let $f$ be a transcendental entire function, and let $V$ be the grand orbit of a component of $\hat F(f)$ contained in a wandering domain. Then, the following are equivalent.
	\begin{enumerate}[label={\normalfont (\alph*)}]
		\item The grand orbit relation in $V$ is discrete.
		\item There exists $z\in V$ such that $\GO(z)$ is discrete in $V$.
		\item For each $z\in V$, there exists a neighbourhood $U$ of $z$ such that $\GO(w)\cap U = \{w\}$ for every $w\in U$.
		\item There exists a simply connected domain $W\subset V$ such that $f|_{f^n(W)}$ is injective for all $n\geq 0$.
	\end{enumerate}
\end{thm}

\section{Dimension of Teichmüller spaces: Proof of Theorem \ref{thm:main}}\label{sec:proofMT}

In this section we provide the proof of our main theorem. Namely, we show that $\T(V,f)$, where $V$  is a component of $\hat F(f)\cap U$ and $U$ a wandering domain of $f$, is infinite dimensional if and only if the orbit relation in $V$ is discrete. It will make use of the concept of \textit{geometric convergence}, both for groups and for marked hyperbolic Riemann surfaces, that we include for convenience (for general references, see \cite{BP92, MT98}).

\begin{defi}[Geometric convergence of groups]\label{def:geolimg}
Let $(\Gamma_n)_{n\in\mathbb{N}}$ be a sequence of subgroups of $\aut(\D)$. 
We say that $(\Gamma_n)$ \emph{converges geometrically} to a subgroup $\Gamma^*< \aut(\D)$ if and only if
	\begin{enumerate}[\rm (a)]
		\item For all $\gamma^*\in\Gamma^*$, there exists $\gamma_n\in \Gamma_n$ such that $\gamma_n\to \gamma^*$ locally uniformly, and
		\item if  $\gamma_{n_k}\in\Gamma_{n_k}$ is a subsequence such that $\gamma_{n_k}\to \gamma\in\aut(\D)$, then $\gamma\in \Gamma^*$.
	\end{enumerate}
\end{defi}

Given a marked hyperbolic Riemann surface $(S,z)$, let $\pi:\D\to S$ be a universal covering map with $\pi(0)=z$. Then $S \simeq\D/\Gamma$, where $\Gamma\le \aut(\D)$ is a discrete, torsion-free subgroup of $\aut(\D)$ known as  the group of {\em deck transformations} of $S$. The choice of $\pi$ is unique up to precomposition by a rotation of $\D$, and correspondingly $\Gamma$ is determined up to conjugation by a rotation in $\aut(\D)$. 

\begin{defi}[Geometric convergence of marked hyp. Riemann surfaces]\label{def:geolimbs}
	Let $((S_n,z_n))_{n\in\N}$  and $(S^*, z^*)$ be marked hyperbolic surfaces and choose universal covers 
	$\pi_n:\D\to S_n$ with $\pi_n(0)=z_n$, and $\pi^*:\D\to S^*$ with $\pi^*(0)=z^*$.
	Let $\Gamma_n,\Gamma^*\le\aut(\D)$ be corresponding groups of deck transformations. We say $((S_n,z_n))$ \emph{converges geometrically} to $(S^*,z^*)$ if  there exist rotations $R_n\in\aut(\D)$ with $R_n(0)=0$ such that 
	$R_n\Gamma_n R_n^{-1}$ converges geometrically to $\Gamma^*$.
\end{defi}


\subsection{Proof of Theorem \ref{thm:main}}
Let $f$ be an entire function with a wandering domain $U$ and let $V$  be a component of $\hat F(f)\cap U$.
If the orbit relation in $V$ is indiscrete, then, by Theorem \ref{thm:MSteich}, $\T(V,f)$ is either trivial or one dimensional, and thus $\T(V,f)$  can only be infinite dimensional if the orbit relation in $V$ is discrete.

For the converse, suppose that the orbit relation in $V$ is discrete. We shall now show that  $\T(V,f)$ is infinite dimensional  by examining the sequence of marked Riemann surfaces $((V_n, z_n))_{n\geq 0}$, where $V_n = f^n(V)$ and $z_n = f^n(z_0)$ for some $z_0\in V$ that we fix. This will be accomplished in three steps; first, we show that this sequence converges geometrically in the sense of Definition \ref{def:geolimbs}.

\begin{lemma}\label{lem:limgeom}
There exists a marked Riemann surface $(V^*, z^*)$ such that $((V_n, z_n))_{n\geq 0}$ converges geometrically to $(V^*, z^*)$.
\end{lemma}
\begin{proof}
For $n\geq 0$, let $\pi_n:\D\to V_n$ be a universal covering map with $\pi_n(0)=z_n$; we allow ourselves to post-compose $\pi_n$ with rotations so as to normalise the argument of $\pi_n'(0)$ when needed. Let $\Gamma_n$ be the group of deck transformations associated to $(V_n,z_n)$. Since $f:V_n\to V_{n+1}$ is a covering, both $f\circ \pi_n$ and $\pi_{n+1}$ are universal covers of $V_{n+1}$. Hence there exists $A_n\in\aut(\D)$ with $\pi_{n+1}\circ A_n = f\circ \pi_n.$ Precomposing $\pi_{n+1}$ with $A_n^{-1}$ (equivalently, adjusting the argument of $\pi_{n+1}'(0)$), and renaming the map, we may assume that the following diagram commutes:
\begin{equation}\label{eq:CD}
	\begin{CD}
		\cdots  (\D,\Gamma_n) @>{\rm Id}>> (\D,\Gamma_{n+1}) \cdots\\
		@V{\pi_n}VV  @VV{\pi_{n+1}}V\\
		\cdots  (V_n,z_n) @>f>> (V_{n+1},z_{n+1}) \cdots
	\end{CD}
\end{equation}

From \eqref{eq:CD} it follows that, for all $n\ge0$, every deck transformation of $\pi_n$ is also a deck transformation of $\pi_{n+1}$; hence $\Gamma_n\subset \Gamma_{n+1}$ and $(\Gamma_n)_n$ is nested. Set
\begin{equation}\label{eq:inclusion}
	\Gamma^*:=\bigcup_{n\ge0} \Gamma_n.
\end{equation}
We shall show that $\Gamma^*$ is a torsion-free discrete subgroup of $\aut(\D)$ and that $(\Gamma_n)_n$ converges geometrically to $\Gamma^*$. We start with the latter. To verify property (a) in Definition~\ref{def:geolimg}, note that for any $\gamma^*\in\Gamma^*$ there exists $N$ such that $\gamma^*\in\Gamma_n$ for all $n\ge N$. Thus the eventually constant sequence $\gamma_n:=\gamma^*\in\Gamma_n$ converges to $\gamma^*$ locally uniformly on compact subsets of $\D$.
To show that property (b) in Definition~\ref{def:geolimg} holds, suppose $\gamma_{n_k}\to\gamma$ with $\gamma_{n_k}\in\Gamma_{n_k}$. We claim that $\gamma\in\Gamma^*$. Consider the points $\gamma_{n_k}(0)\in\D$. Since $\pi_{n_k}\circ\gamma_{n_k}=\pi_{n_k}$ and $\pi_{n_k}(0)=z_{n_k}$, we have $\pi_{n_k}(\gamma_{n_k}(0))=z_{n_k}$. For each $k$, define
\[
w_k:=\pi_0(\gamma_{n_k}(0))\in V_0.
\]
Then, using the commutative diagram \eqref{eq:CD} (so that $\pi_{n_k}=f^{n_k}\circ\pi_0$), we see that   $f^{n_k}(w_k)=z_{n_k}$, so $w_k$ lies in the grand orbit of $z_0$. Because $\gamma_{n_k}\to\gamma$ locally uniformly, $\gamma_{n_k}(0)\to v^*:=\gamma(0)$ in $\D$, hence $w_k\to w^*:=\pi_0(v^*)$ in $V_0$. Since the grand orbit relation in $V_0$ is discrete, by Theorem~\ref{thm:GOrelations}. any convergent sequence in $\GO(z)$ must be  eventually constant;  thus there exists $K$  such that $w_k=w^*$ for all $k\ge K$.

Fix such $k\ge K$. Then $\pi_0(\gamma_{n_k}(0))=\pi_0(v^*)=w^*$. As $\pi_0$ is a covering and its fiber over $w^*$ is discrete, the convergence $\gamma_{n_k}(0)\to v^*$ forces $\gamma_{n_k}(0)=v^*$ for all sufficiently large $k$. Passing to this tail, we have
\[
\gamma_{n_k}(0)=v^*\quad\text{and}\quad \gamma_{n_k}\to\gamma.
\]
It follows that the tail is actually constant; indeed, $\gamma_{n_k}\circ \gamma^{-1}_{n_{k+1}}$ is a rotation and belongs to $\Gamma_{n_{k+1}}$, which is not possible since this group is discrete and torsion-free. In particular, $\gamma_{n_k}=\gamma$ for all $k$ large enough. Since each $\gamma_{n_k}\in\Gamma_{n_k}\subset\Gamma^*$, we conclude $\gamma\in\Gamma^*$.

It remains to prove that $\Gamma^*$ is torsion-free and discrete, ensuring that $V^*:= \D/\Gamma^*$ is a hyperbolic Riemann surface. Since each $\Gamma_n$ is torsion-free, so is their union $\Gamma^*$. For discreteness, suppose not. Then, there exists a sequence $\gamma_k\in\Gamma^*$ with $\gamma_k\neq\Id$ and $\gamma_k\to\Id$ locally uniformly on $\D$. Arguing as above, let $n(k)$ be minimal with $\gamma_k\in\Gamma_{n(k)}$. Set $u_k:=\gamma_k(0)\in\D$ and $w_k:=\pi_0(u_k)\in V_0$. Since $\gamma_k\to\Id$, we have $u_k\to 0$, hence $w_k\to z_0=\pi_0(0)$. Because $\pi_{n(k)}=f^{\,n(k)}\circ\pi_0$ and $\pi_{n(k)}\circ\gamma_k=\pi_{n(k)}$, we obtain
\[
f^{\,n(k)}(w_k)=\pi_{n(k)}(u_k)=\pi_{n(k)}(\gamma_k(0))=\pi_{n(k)}(0)=z_{n(k)}.
\]
Thus each $w_k$ lies in the grand orbit of $z_0$. Since the grand orbit relation in $V_0$ is discrete, the convergence $w_k\to z_0$ forces $w_k=z_0$ for all large $k$. As the fiber $\pi_0^{-1}(z_0)$ is discrete and $u_k\to0$, it follows that $u_k=0$ for all large $k$, i.e. $\gamma_k(0)=0$ eventually. But a nontrivial automorphism of $\D$ fixing $0$ is a rotation. If its angle is rational, it is finite order and is excluded by torsion-freeness; if its angle is irrational, then its powers accumulate at the identity, contradicting discreteness. Hence no nontrivial sequence in $\Gamma^*$ can converge to $\Id$, and $\Gamma^*$ is discrete.

We have thus shown that $(\Gamma_n)_n$ converges geometrically to $\Gamma^*$. By Definition~\ref{def:geolimbs}, it follows that the marked hyperbolic Riemann surfaces $((V_n,z_n))_n$ converge geometrically to the marked surface $(V^*,z^*)$, where $V^*=\D/\Gamma^*$ and $z^*:=\pi^*(0)$ for the quotient map $\pi^*:\D\to V^*$.
\end{proof}

Next, we must relate this limit to the dynamics and, most importantly, the Teichm\"uller space of $f$.

\begin{lemma}\label{lem:isom_surfaces}
The surface $V^*$ is isomorphic to $\GO(V)/f$, and hence $\T(V, f) \simeq \T(V^*)$.
\end{lemma}
\begin{proof}
	Observe that $\GO(V)/f$ is a Riemann surface because the grand orbit relation in $V$ is discrete (Theorem \ref{thm:MSteich}). Also, for every $n\ge0$ we have the projections\footnote{There is a slight abuse of notation here: formally, the operation ``$/f$'' denotes a covering map $\Psi\colon \bigsqcup_n \hat V_n\to \GO(V)/f$, where $\hat V_n$ runs over all connected components of $\GO(V)$, and we are taking its restriction to $V_n$.} 
	\[
	(\D,\Gamma_n)\xrightarrow{\ \pi_n\ }V_n \xrightarrow{/f} \GO(V)/f,
	\]
	and, by the commutative diagram \eqref{eq:CD}, $\pi_n=f^n\circ\pi_0$.
	
	To see that $\D/\Gamma^*\simeq \GO(V)/f$, first let $[z]=[w]\in\D/\Gamma^*$. Then $\gamma(z)=w$ for some $\gamma\in\Gamma^*$. Since $\gamma\in\Gamma_n$ for some $n$ by \eqref{eq:inclusion}, we have $\pi_n(z)=\pi_n(w)\in V_n$; hence $z$ and $w$ project trivially to the same point in $\GO(V)/f$. Conversely, suppose $z,w\in\D$ project to the same point in $\GO(V)/f$. Then there exists $n\ge0$ with $f^{\,n}(\pi_0(z)) = f^{\,n}(\pi_0(w)).$ By \eqref{eq:CD}, $\pi_n=f^{\,n}\circ\pi_0$, hence
	$\pi_n(z)=\pi_n(w)$,
	so there exists $\gamma\in\Gamma_n\subset\Gamma^*$ with $\gamma(z)=w$.
	
	We have shown the equivalence relation induced by $\Gamma^*$ on $\D$ agrees with the grand-orbit quotient under $f$, so the induced map $\D/\Gamma^*\to \GO(V)/f$ is a bijection. As both sides are Riemann surfaces and the map is induced by holomorphic coverings, it is a biholomorphism. 
	
	Finally, by Theorem~\ref{thm:MSteich} (applied to $\widehat V=\GO(V)$), the discreteness of the grand orbit relation yields $\T(V,f)\cong \T(\GO(V)/f)=\T(V^*),$ as we wanted to show.
\end{proof}

Finally, recall that a Riemann surface has a finite-dimensional Teichm\"uller space if and only if it is of finite type \cite[Proposition 7.1.1]{Hub06}. Since each $V_n$ is of infinite type (it is an open set with non-isolated boundary points),  discreteness of the orbit relation in $V$ implies that $\T(V, f)$ is infinite dimensional; this follows by combining Lemma \ref{lem:isom_surfaces} with the next one.
\begin{lemma}[Geometric limit of surfaces of infinite type] \label{lem:inflimits}
Let $f_n\colon (S_n, p_n)\to (S_{n+1}, p_{n+1})$ be a sequence of holomorphic coverings between hyperbolic Riemann surfaces. Assume that the sequence $((S_n, p_n))_{n\geq 0}$ converges geometrically to $(S^*, p^*)$, and that each $S_n$ is of infinite type. Then $S^*$ is of infinite type.
\end{lemma}

\begin{remark}
A version of this lemma, with a different terminology, can be found in \cite[Proposition 1]{sullivan}. We give here a short proof for completeness.
\end{remark}
\begin{proof}
	Since $f_n\colon (S_n,p_n)\to (S_{n+1},p_{n+1})$ is a covering, the induced map $(f_n)_*\colon \pi_1(S_n)\mapsto \pi_1(S_{n+1})$ is injective. Equivalently, identifying each $\pi_1(S_k)$ with the group of deck transformations of the universal cover of $S_k$, we have (after fixing basepoints and the usual identifications) that $\pi_1(S_n)\;\subset\;\pi_1(S_{n+1})$ for all $n$.
	 Furthermore, it holds that 
\[ \pi_1(S^*)\simeq \bigcup_{n\geq 1} \pi_1(S_n),\]
see the proof of Lemma \ref{lem:limgeom} for details. Now, if $S^*$ is of finite type, then $\pi_1(S^*)$ is finitely generated; for $n$ large enough, the fundamental group of $S_n$ contains all the generators of $\pi_1(S^*)$, and so, since $\pi_1(S_n)\subset \pi_1(S^*)$, $\pi_1(S_n)$ is also finitely generated -- a contradiction.
\end{proof}
This concludes the proof of Theorem \ref{thm:main}.

\section{Linearising coordinates: proof of Theorem \ref{thm:discrete_intro}}\label{sec:proofdiscrete}
In this section we provide the proof of Theorem \ref{thm:discrete_intro}, that is, we provide a characterisation of discrete orbit relation of a component  of $\hat{F}$ in a wandering domain in terms of the existence of what can be thought as linearising coordinates. Since the orbit relation in any component of a multiply connected wandering domain is indiscrete, see \cite[Corollary 6.9]{EFGP24}, our discussion in this section will only concern simply-connected wandering domains. We will make use of this observation without further mentioning.

\subsection{Complex analytical tools: distortion estimates}
We shall require the following collection of results and estimates of holomorphic maps, that we include for the reader's convenience.  
\begin{lemma}[Zalcman's Lemma \cite{Zal75}]\label{lem:zalcman}
	A family $\mathcal{F}$ of functions holomorphic in $\DD$ is \emph{not} normal in any neighbourhood of $z_0\in\DD$ if and only if there exist:
	\begin{enumerate}[(a)]
		\item A sequence $(z_n)_n\subset\DD$ such that $z_n\to z_0$;
		\item A sequence $(\rho_n)_n$ of positive real numbers such that $\rho_n\to 0$; and
		\item A sequence of functions $(f_n)_n\subset \mathcal{F}$ such that
		\[ f_n(z_n + \rho_nz) \underset{n\to\infty}{\longrightarrow} f(z) \]
		locally uniformly in $\Cx$, where $f\colon\Cx\to\Cx$ is a \emph{non-constant} entire function.
	\end{enumerate}
\end{lemma}
We will make use of the following version of Koebe's distortion theorem:
\begin{lemma}[Koebe Distortion Theorem]\label{lem:koebe}
	Let $f\colon\DD_R\to\Cx$ be univalent with $f(0) = 0$. Then, for any $z\in\DD_R$,
	\[ \frac{R^2|z||f'(0)|}{(R + |z|)^2} \leq |f(z)| \leq \frac{R^2|z||f'(0)|}{(R - |z|)^2} \]
\end{lemma}

Finally, we will need  the following estimates on self-maps of the unit disc.
\begin{lemma}[\cite{BEFRS21}, Corollary 2.4]\label{lem:BEFRS}
	Let $f\colon\DD\to\DD$ be holomorphic with $f(0) = 0$ and $|f'(0)| = \lambda$. Then, for any $z\in\DD$,
	\[ |z|\left(\frac{\lambda - |z|}{1 - \lambda|z|}\right) \leq |f(z)|. \]
\end{lemma}

\begin{lemma}\label{lem:univ}
For every $c > 0$, there exists $R > 0$ such that every holomorphic map $g\colon\DD\to\DD$ with $g(0) = 0$ and $|g'(0)| > 0$ is univalent on $\DD_R$.
\end{lemma}
\begin{proof}
	First note that for any critical point $w$ of $g$,
	\[ 
	2d_\DD(0, w) \geq d_\DD(\|Dg(0)\|_\DD^\DD, \|Dg(w)\|_\DD^\DD) = d_\DD(g'(0), 0) \geq d_\DD(c, 0), 
	\]
	where the first inequality is Schwarz--Pick inequality for hyperbolic distortions (\cite[Theorem 11.2]{BM06}). This means that there exists $R_1 =R_1(c)> 0$ such that any critical point $w$ of $g$ satisfies $|w| \geq R_1$. Next, define
	\[
	h_c(r):=r\,\frac{c-r}{1-cr}\qquad r\in[0,1),~c\in(0,1].
	\]
	By Lemma~\ref{lem:BEFRS} we have \(|g(z)|\ge h_{|g'(0)|}(|z|)\) for all \(z\in\mathbb{D}\). Because \(c\le |g'(0)|\) and, for fixed \(r\), the function \(c\mapsto h_c(r)\) is strictly increasing, it follows that \(|g(z)|\ge h_c(|z|)\) for all \(z\in\mathbb{D}\). For any $r\in (0,1)$, by Schwarz's lemma, $g^{-1}(\DD_r)$ has a connected component $D$ such that $D\supset \DD_r$. 	Moreover, if $0<r<h_c(R_1)$, then for every $z$ with $|z|=R_1$ we have $|g(z)|\ge h_c(R_1)>r$, so the circle $\{|z|=R_1\}$ does not meet $g^{-1}(\mathbb D_r)$. Hence $D\subset\mathbb D_{R_1}$. Since $g$ has no critical points in $\mathbb D_{R_1}$, the map $g:D\to\mathbb D_r$ is an unbranched covering; as $\mathbb D_r$ is simply connected, $g$ is a biholomorphism on $D$, and in particular univalent on $\mathbb D_r$. Therefore $g$ is univalent on $\mathbb D_r$ for every $0<r<h_c(R_1)=:R$, where $R$ depends only on~$c$.		
\end{proof}

Finally, we will also use the following fact on convergence of real sequences.
\begin{obs}\label{obs:cauchy}\normalfont Let $(x_n)_{n\in\mathbb{N}}$ be a sequence of positive numbers such that $x_n \leq M$ for all $n\in\mathbb{N}$. Then, if $\frac{x_{m+k}}{x_m}\to 1$ as $m\to+\infty$ uniformly in $k\in\mathbb{N}$, the sequence $(x_n)_{n\in\mathbb{N}}$ converges. Indeed, for any $\epsilon>0$ choose $N$ so that $\left|1-\frac{x_{m+k}}{x_m}\right|<\epsilon/M$ for all $m\ge N$ and all $k\in\mathbb{N}$. Then
\[
|x_{m+k}-x_m|=\left|x_m\left(1-\frac{x_{m+k}}{x_m}\right)\right|
\le M\left|1-\frac{x_{m+k}}{x_m}\right|<\epsilon,
\]
so $(x_n)$ is Cauchy and hence convergent.
\end{obs}

\subsection{Proof of Theorem \ref{thm:discrete_intro}}
Let $f$ be an entire function with a wandering domain $U$,  and for each $n\in \N$, let $U_n$ be the Fatou component containing $f^n(U)$. Moreover, let $V=:V_0$ be a component of $\widehat F(f)\cap U$ and denote $V_n:=f^n(V)$ for all $n\in \N$. By Theorem~\ref{thm:GOrelations}, the orbit relation in $V$ being discrete (item \ref{item:discrete} in the statement of Theorem~\ref{thm:discrete_intro}) is equivalent to the following.
	\begin{enumerate}[(a')]
		\item \label{item:c_proof}  For each $z\in V_0$ and all $n\geq 0$ there exists a neighbourhood $W_n\subset V_n$ of $f^n(z)$ such that $f$ maps $W_n$ to $W_{n+1}$ injectively. 
	\end{enumerate}
Thus, we will equivalently prove that \ref{item:linearising} $\Leftrightarrow$ \ref{item:c_proof}. Observe that the implication \ref{item:linearising} $\Rightarrow$ \ref{item:c_proof}  is direct, since  the assumptions in \ref{item:linearising} imply that $f|_{W_n}$ is injective for all $n$.
	
For the converse, let us assume that \ref{item:c_proof} holds (which also means that the orbit relation in $V$ is discrete), let us fix a point $\tilde{z}_0\in V$, denote $\tilde{z}_n=f^n(\tilde{z}_0)\in V_n$ for each $n$, and let $(W_n)_n$ be the corresponding neighbourhoods of injectivity given by \ref{item:c_proof}. For short, we will denote $\lambda_n=\lambda_n(\tilde{z}_n)= \|Df(f^{n}(z_0))\|_{U_n}^{U_{n+1}}.$ Following subsection \ref{subsec:wd}, let $(g_n)$ be a sequence of inner functions associated with $f$ in $(U_n)$ via the Riemann maps 
$$\psi_n\colon U_n\to \D^{(n)},\quad \text{with } \psi_n(\tilde{z}_n)=0 \text{ for all } n,$$
and recall from \eqref{eq_deriv0} that $\la_n=\vert g_n'(0)| \neq 0$ for all $n\geq 0$, since $\tilde{z}_n\in \widehat F(f)$.

	Given $m,n\geq 0$ such that $m > n$, define $H_n^m\colon V_n\to \C$ as
	\[ 
	H_n^m(z) := \frac{\psi_m(f^{m-n}(z))}{\Lambda_n^m}, \quad
	\text{ where } \Lambda_n^m := \lambda_{m-1}\cdots\lambda_n.
	\] 
	Observe that, by definition, the map $H_n^m$ is holomorphic and satisfies $H_n^m(\tilde{z}_n)=0$. 
	
	Next, we aim to show that for each $n$, the sequence $(H_n^m)_{m\geq 1}$ converges locally uniformly as $m\to\infty$ by showing that the sequence forms a normal family and that there is a unique limit function, $\phi_n$,  such that  the sequence $(H_n^m)_m$ converges locally uniformly to $\phi_n$. Additionally, we will show that $(\phi_n)$ satisfies the properties stated in \ref{item:linearising}.
	
	\begin{lemma}\label{lem:normal} For each $n\in \N$, the family $(H^m_n)_m$ is normal in $V_n$.
	\end{lemma}
	\begin{proof}
		Let us assume, for the sake of contradiction, that $(H^m_n)_m$ is not normal in any neighbourhood of some $z_0\in V_n$. Then, by Zalcman's lemma (Lemma \ref{lem:zalcman}), there exist $z_k\to z_0$, $\rho_k\to 0$, and $(m_k)_k$ such that the sequence $(h_k)_k$ given by $h_k(w) := H_n^{m_k}(z_k + \rho_kw)$ converges to a non-constant entire function $h$. We consider two cases, according to whether $h$ is affine or not.
		
		If $h$ is affine then it has a zero $w^*$, and by Hurwitz's theorem there exists a sequence $w_k\to w^*$ such that $h_k(w_k) = 0$ for all $k$. If we define $p_k := z_k + \rho_kw_k$, we see that $p_k\to z_0$ as $k\to\infty$, and also that $H_n^{m_k}(p_k) = 0$ for all $k$. Then, by definition of the maps involved, 
		$$H_n^{m_k}(p_k)=0 \;\Rightarrow\; \psi_{m_k}\!\big(f^{\,m_k-n}(p_k)\big)=0 \;\Rightarrow\; f^{\,m_k-n}(p_k)=\tilde z_{m_k}=f^{\,m_k-n}(\tilde z_n) \text{ for all }k.$$
		In particular, $p_k \in \GO(\tilde{z}_n)$ for all $k$, which implies that any neighbourhood of $z_0$ has infinitely many points of $\GO(\tilde{z}_n)$. This contradicts Theorem~\ref{thm:GOrelations}(c).
		
		Similarly, if $h$ is not affine, then it has at least one singular value $s\in\Cx$; by continuity of the singular values \cite[Lemma 1]{KK97}, there exists a sequence $s_k\in S(h_k)$ such that $s_k\to s$ as $k\to\infty$. If we define $p_k := z_k + \rho_ks_k$, we see that once again $p_k\to z_0$ as $k\to\infty$, and that each $p_k$ is a singular point (mapped to a singular value) of $H^{m_k}_n$, and hence of the composition $f^{m_k-n}$. This would imply that $z_0\in \overline{\GO(S(f))}$, which contradicts that $z_0\in V_n\subset \hat{F}(f)$.
	\end{proof}
	
	\begin{lemma}\label{lem:pointwise}
		For each $n\in \N$ there exists $R=R(n)>0$ such that the sequence $(|H_n^m(z)|)_{m\geq n}$ converges pointwise for each $z\in \D_R(\tilde{z}_n)$.
	\end{lemma}
	
	\begin{proof}
		In order to see this, we will show the analogous property for the maps $E^m_n:= H^m_n\circ \psi_n^{-1}\colon X_n\subset \D^{(n)}\to \C$, where $X_n=\psi_n(V_n)$. That is, we will show that there exists $Q>0$ such that $(|E_n^m(z)|)_{m\geq n}$ converges pointwise for each $z\in \D_Q$. Then, taking $R>0$ such that $\psi_n(\D_R(\tilde{z}_n))\subset \D_Q$, the result will follow. Note that by definition, 
		\begin{equation}\label{eq:Enm}
				E_n^m(z)=\frac{G_n^m(z)}{\Lambda_n^m}, \quad
		\text{ where } G_n^m(z) := g_{m-1}\circ\cdots\circ g_n(z) \in \DD^{(m)}.
		\end{equation}
		Observe that for all $m,n$ with $m\geq n$, the maps $E_n^m$ are holomorphic, fix the origin and satisfy $\vert (E_n^m)'(0)\vert= 1$.
		Since we have assumed that the orbit relation in $V$ is discrete, $U$ is simply-connected and so the classification in Theorem \ref{thm:Ainternal} applies. We divide rest of the proof into two cases, depending on whether $U$ is contracting or not. 
		
		If $U$ is not contracting, then it is either eventually isometric or semi-contracting. In these cases, by Theorem \ref{thm:Ainternal},
the distortions $\lambda_n \in (0,1]$ tend to one as $n\to \infty$. Moreover, as $m\to \infty$ the constants $\Lambda_n^m$ tend to a non-zero limit 
		\[ 
		\Lambda_n := \prod_{m\geq n} \lambda_m. 
		\] 
Indeed, this follows from standard estimates of infinite products (see e.g. \cite[Theorem~15.5]{Rud87}), which establish that if $\lambda_n\in (0,1]$ 
	  \begin{equation}\label{eq:inf_prod}
	  	\sum_{n=0}^{\infty} (1 - \lambda_n) < \infty \iff \prod_{n=0}^{\infty} \lambda_n>0,
	  \end{equation}
		and the equality in the left hand side, by Theorem \ref{thm:Ainternal}, is equivalent to $U$ not being contracting.
		
		Since each $g_n$ maps the disc onto itself, it follows that $|E_n^m(z)| \leq (\Lambda_n)^{-1}$ for all $z\in X_n$ and $m\in\mathbb{N}$. Choose $Q\leq \vert \Lambda_1 \vert$ such that $\D_Q\subset X_n$. Fix $k\in\mathbb{N}$, and observe that $\Lambda_{m+1}^{m+k}\to 1$ as $m\to+\infty$, since
		 		\[\Lambda_{m+1} \leq  \Lambda_{m+1}^{m+k}\leq 1, \]
		  and $\lambda_n$, and so  $\Lambda_{m+1}$, tends to $1$ as $m\to+\infty$.
		Now, by Lemma~\ref{lem:BEFRS} applied to the map $G^{m+k}_{m+1}$ and by Schwarz's lemma, we have
		\begin{equation}\label{eq:semi}
			|z|\left(\frac{\Lambda_{m+1}^{m+k} - |z|}{1 - \Lambda_{m+1}^{m+k}\cdot |z|}\right) \leq |G_{m+1}^{m+k}(z)| \leq |z|.
		\end{equation}
		If $z\in\DD_Q$, then $\Lambda_{m+1}^{m+k} - |z|>Q-|z|>0$ and hence, taking the limit $m\to+\infty$ we see that $|G_{m+1}^{m+k}(z)|\to |z|$ for every $z\in\DD_{Q}$. It also follows from \eqref{eq:semi} that $G_{m+1}^{m+k}(z)$ is non-zero for $z\in\DD_{Q}\setminus\{0\}$ for any $m$ and $k$. Now, by definition of the map $E_n^{m+k}$,
		\begin{equation}\label{eq:E_quotient}
\frac{|E_n^{m+k}(z)|}{|E_n^m(z)|} = \frac{|G_n^{m+k}(z)|}{|G_n^m(z)|}\frac{1}{\Lambda_{m+1}^{m+k}} = \frac{|G_{m+1}^{m+k}(G_n^m(z))|}{|G_n^m(z)|}\frac{1}{\Lambda_{m+1}^{m+k}}.
		\end{equation}
		Since the right-hand side of the equation goes to one uniformly on $k$ as $m\to+\infty$, by (\ref{eq:semi}) and since $\Lambda_{m+1}^{m+k}\to 1$, we reach the desired conclusion using Observation~\ref{obs:cauchy}.

		We are left to study the case when $U$ is contracting. Assume so and let $Q> 0$ such that $\D_Q\subset \psi_n(W_n)$. In particular, since $G^m_n=\psi_m\circ f^{m-n}\circ \psi^{-1}_n$, the maps $G^m_n$ and $E^m_n$ are univalent in $\D_Q$ for all $m$. Let us fix $k\in \N$. Applying Koebe distortion theorem (Lemma~\ref{lem:koebe}) to the map $G^{m+k}_{m+1}$ restricted to $G_n^m(\DD_Q) \subset \DD_Q$ and using the expression obtained in \eqref{eq:E_quotient}, it holds that
		\begin{equation}\label{eq:appl_Koebe}
		\frac{Q^2}{(Q+ |G_n^m(z)|)^2} \leq \frac{|E_n^{m+k}(z)|}{|E_n^m(z)|} \leq \frac{Q^2}{(Q- |G_n^m(z)|)^2} \quad \text{ for all }z\in \D_Q.
		\end{equation}
		By the characterisation of $U$ being contracting in terms of divergence of the series given in Theorem \ref{thm:Ainternal} combined with \eqref{eq:Gnto0}, $G_n^m(z)\to 0$ as $m\to+\infty$. Combining this with \eqref{eq:appl_Koebe}, we can apply Observation~\ref{obs:cauchy} and conclude that $|E_n^m(z)|$ converges as $m\to+\infty$ for any $z\in\DD_Q$. 
	\end{proof}

	It follows from Lemmas \ref{lem:normal} and \ref{lem:pointwise} that for each $n$, any two limit functions of the sequence of maps $(H^m_n)_m$ restricted to $\D_{R(n)}(\tilde{z}_n)$ differ in at most a rotation. Moreover, 
	$$(H_n^m)'(\tilde{z}_n)=\frac{(\psi_m(f^{m-n}(\tilde{z}_n)))'}{\Lambda_n^m}= \frac{(G^m_n(\psi_n(\tilde{z}_n)))'}{\Lambda_n^m}=\frac{\Lambda_n^m\cdot\psi'_n(\tilde{z}_n)}{\Lambda_n^m}=\psi'_n(\tilde{z}_n),$$
	where we have used $m$ times the commutative relation \eqref{eq:associated_inner} and that $(G^m_n(\psi_n(\tilde{z}_n))'=(G^m_n(0))'=\Lambda_n^m$ by \eqref{eq_deriv0}, $G^m_n$ being defined in \eqref{eq:Enm}.
	 In particular, the value of $(H_n^m)'(\tilde{z}_n)$ does not depend on $m$, and so all limit functions must coincide. In other words, we have shown that 
	\[ 
	\phi_n(z) :=  \lim_{m\to+\infty} H_n^m(z) 
	\]
	exists when restricted to $\D_{R(n)}(\tilde{z}_n)$. Now, since we saw in Lemma \ref{lem:normal} that $(H_n^m)$ is normal in $V_n$ it follows from the partial converse of Montel's Theorem (see e.g. \cite[Corollary 2.2.4]{schiff}) that it is locally bounded, i.e. uniformly bounded in a neighborhood of each point of $V_n$. Hence, by Vitali convergence theorem (see e.g. \cite[Section 2.4]{schiff}), the sequence converges locally uniformly on all $V_n$, and so the definition of $\phi_n$ extends to all $z\in V_n$. 

	Finally, observe that by construction, for each $z\in V_n$,
	\[ \phi_{n+1}\circ f(z) = \lim_{m\to+\infty} \frac{\psi_m(f^{m-(n+1)}(f(z))}{\Lambda_{n+1}^m} = \lambda_n\cdot\lim_{m\to+\infty} \frac{\psi_m(f^{m-(n+1)}(f(z))}{\Lambda_n^m} = \lambda_n\cdot\phi_n(z), \]
	as desired. Since for each $m$, the restriction of $H^m_n$ to $W_n$ is univalent as composition of univalent functions and the limit function $\phi_n$ is not constant, it follows that $\phi_n$ is univalent in $W_n$, which concludes the proof of the theorem. 

\begin{remark}We emphasise again that the maps $\phi_n$ can be post-composed with arbitrary linear maps to produce new linearising coordinates, conjugating $f|_{\textcolor{red}{V}_n}$ to a different sequence of linear maps.
\end{remark}

\section{Teichmüller spaces: Proof of Theorem \ref{thm:teich}}\label{sec:proofteich}

In the proof of Theorem \ref{thm:discrete_intro}, many of our arguments concerned self-maps of the disc rather than wandering domains. These arguments can be adapted to obtain an analogous result for sequences $(g_n)$ of inner functions fixing the origin whose derivative is uniformly bounded away from zero, which can be seen again as a non-autonomous version of the usual K\"onigs coordinates for the ``fixed point'' at the origin. 
 We remark that this result was already proven by Pommerenke in \cite{Pom94}, with a different proof. Since this theorem will be relevant in the proof of Theorem \ref{thm:teich}, we outline here our version of the proof based on that of Theorem \ref{thm:discrete_intro}, indicating where the arguments differ and what additional considerations must be taken into account to adapt the proof to this setting.

\begin{thm}[Uniform linearising coordinates for inner functions]\label{thm:simlin}
	Let $g_n\colon\DD\to\DD$, $n\geq 0$, be inner functions  fixing the origin, and let $\lambda_n:= |g_n'(0)|$. If there is a constant $c > 0$ such that $\lambda_n \geq c$ for all $n\in\mathbb{N}$, then there exists a sequence $\phi_n\colon\DD\to\Cx$ of analytic functions fixing the origin such that 
	\[ \phi_{n+1}\circ g_n(z) = \lambda_n\cdot\phi_n(z) \quad \text{ for all }n\geq 0 \text{ and } z\in \D. \]
	Furthermore, there exists $R > 0$ such that $\phi_n|_{\DD_R}$ is univalent for all $n\in\mathbb{N}$, and the sequence $(\phi_n)_n$ is uniquely determined by the additional requirement that~$\phi_n'(0) = 1$.
\end{thm}
\begin{proof}[Outline of the proof] Firstly we note that, by Lemma \ref{lem:univ}, the moduli of the derivatives at zero being uniformly bounded implies that there exists a uniform radius $R>0$  where $g_n|_{\DD_R}$ is univalent for all $n$. For each $n$ and $m\geq n$ and $z\in \D$, define
	\begin{equation} \label{eq:E2}
 E_n^m(z)=\frac{G_n^m(z)}{\Lambda_n^m}, \quad
	\text{ where } G_n^m(z) := g_{m-1}\circ\cdots\circ g_n(z) \text{ and }\Lambda_n^m := \lambda_{m-1}\cdots\lambda_n.
	\end{equation}
	By analogy with the wandering domains case and using \eqref{eq:Gnto0}, we will say that $(G^m_n)_{m}$ is \textit{contracting} if $\Lambda_n=:\prod_{m\geq n} \lambda_n$ converges to zero (see \cite[Theorem 7.2]{BEFRS22}).  We can divide the proof into two cases, depending on whether $(G^m_n)_{m}$ is contracting or not. 
	
	Arguing as in the proof of Lemma \ref{lem:pointwise}, one can see that $(|E_n^m(z)|)_{m\geq n}$ converges pointwise for each $z\in \D_Q$, where $Q=\vert \Lambda_1\vert $ if $(G^m_n)_{m}$ is not contracting, and $Q=R$ if $(G^m_n)_{m}$ is contracting. 
	
	Next, if $(G^m_n)_{m}$ is not contracting, then, since $|E_n^m(z)\vert<(\Lambda_n)^{-1}$ for all $z\in \D$ and $m\in \N$, we have by Montel's theorem that $(E_n^m)_m$ is a normal family in $\D$. However, in the contracting case,  we can only conclude normality of $(E_n^m)_m$ in $\D_R$ by using Koebe distortion theorem (Lemma \ref{lem:koebe}) to show that for $z\in \D_R$,
	\begin{equation}
		\vert E_n^m(z)\vert \leq\frac{R^2\vert z \vert}{(R-\vert z \vert)^2},
	\end{equation}
	together with Montel's theorem.
	
	Arguing as in the proof of Theorem \ref{thm:discrete_intro}, we conclude that 
	\begin{equation}\label{eq:lim5}
\phi_n(z) :=  \lim_{m\to+\infty} E_n^m(z)
	\end{equation} 
	is well defined for all $z\in\D$ in the non-contracting case, and for all $z\in \D_R$ in the contracting case. To finish the proof in the contracting case, we can extend the definition of $\phi_n$, $n\in \N$, to all $z\in \D$ by \textit{spreading by the dynamics}: in this case, for all $z\in \D^{(n)}$, there exists $m$ such that $G_n^m(z)\in \D_R$, and thus, we can define  $\phi_n(z)=\frac{\phi_{n+m}(G_n^m(z))}{\Lambda_n^m}.$
\end{proof}

\subsection{Proof of Theorem \ref{thm:teich}} 
Let $f$ be an entire function with a wandering domain $U$, for each $n\in \N$, let $U_n$ be the Fatou component containing $f^n(U)$ and call $U_0=U$. Following the hypotheses of the statement, there exists
 $z\in U$ for which the hyperbolic distortions $\lambda_n = \|Df(f^{n}(z))\|_{U_{n}}^{U_{n+1}}$ satisfy $\lambda_n \geq c > 0$.\footnote{Combining Theorem~\ref{thm:Ainternal} with the Schwarz--Pick inequality for hyperbolic distortions \cite[Thm.~11.2]{BM06}, we conclude that the same property holds for every $w\in U$.}
As described in subsection \ref{subsec:wd}, 
taking Riemann maps $\psi_n\colon U_n\to \D$ such that $\psi_n(f^n(z_0))=0$, we can define a sequence of inner functions $g_n\colon \D\to\D$ satisfying 
\begin{equation*}\label{eq:associated_inner2}
\psi_{n+1} \circ f=g_n \circ \psi_n,
\end{equation*}
fixing the origin and so that $g_n'(0) = \lambda_n$. Then, by Theorem \ref{thm:simlin}, there exist holomorphic maps $\phi_n\colon \DD\to\Cx$ such that $$\phi_{n+1}\circ g_n=\lambda_{n}\cdot \phi_n,$$ and these maps are all univalent in a disc of radius $r > 0$ around $0$. Hence, the compositions $\varphi_n := \phi_n\circ \psi_n\colon U_n\to \D$, $n\geq 0$, satisfy 
\begin{equation}\label{eq:comp}
\varphi_{n+1}\circ f= \lambda_n\cdot \varphi_n \quad \text{ and } \quad \varphi_n(f^n(z_0))=0.
\end{equation}

Let $W_n:=\psi_n^{-1}(\D_r)$. Since $\phi_n$ is univalent on $\D_r$, the chart $\varphi_n=\phi_n\circ\psi_n$ is univalent on $W_n$. Moreover, by Schwarz–-Pick, $g_n(\D_r)\subset\D_r$, hence $f(W_n)\subset W_{n+1}$. Thus, we have found \textit{linearising coordinates} for the dynamics of $(f|_{U_n})$ in the sense of Theorem \ref{thm:discrete_intro}, and so the orbit relation in the whole of $U$ is discrete\footnote{This does not imply $U\cap \hat F(f)$ is connected: preimages of $U$ can introduce points of the {\em postsingular set} $P(f)$ (forward images of $S(f)$) inside $U$, and the set $(P(f)\setminus S(f))\cap U$ may carve $U$ into several components.}

\begin{lemma} \label{lem:GOphi}
For $z, w\in U_n$ and $n\geq 0 $, $\varphi_n(z) = \varphi_n(w)$ if and only if $\GO(z) = \GO(w)$. 
\end{lemma}
\begin{proof}
It follows directly from the commutative relation \eqref{eq:comp} that if $z$ and $w$ are in the same grand orbit, that is, $f^m(z)=f^m(w)$ for some $m \geq 1$, then $\varphi_n(w) = \varphi_n(z)$. For the converse, assume that $z$ and $w$, both in $U_n$,  satisfy $\varphi_n(z) = \varphi_n(w)$. Let $p=f^n(z_0)\in U_n$, and for each $m\geq 0$, let $r_m > 0$ be the radius of the largest (hyperbolic) disc around $f^m(p)$ such that $\varphi_{n+m}\vert_{D_{U_{n+m}}(f^m(p), r_m)}$ is injective, which, in particular, is contained in $W_n$. Notice that Theorem \ref{thm:simlin} gives us a positive lower bound for $r_m$, and that, furthermore, if $U$ is not contracting, then $r_m\to\infty$ as $m\to\infty$.

Suppose for the sake of contradiction that $f^m(z)\neq f^m(w)$ for all $m$. By the Schwarz--Pick inequality $d_{U_{n+m}}(f^m(p),f^m(z))$ is decreasing in $m$, and the same is true for $w$. In the contracting case it decreases to 0 and hence the lower bound in $r_m$ ensures that for all  $m\in\mathbb{N}$ sufficiently large, the points $f^{n+m}(z)$ and $f^{n+m}(w)$ both belong to $D_{U_{n+m}}(f^m(p), r_m)$. The same is true in the non-contracting case, because $r_m\to\infty$. But  since $\varphi_{n+m}$ is injective in this disc, the points $\varphi_{n+m}(z)$ and $\varphi_{n+m}(w)$ must also be different, which, in turn, by the commutative relation \eqref{eq:comp}, implies that   $\varphi_n(z)$ and $\varphi_n(w)$ are also different points -- a contradiction.
\end{proof}

By Lemma~\ref{lem:GOphi} for $n=0$, the map \(\varphi_0\) is constant precisely on grand orbits and separates distinct ones; hence it induces a biholomorphism
\[
(\GO(U)\cap \hat F(f))/f \;\xrightarrow{\;\sim\;}\; \varphi_0(U)\setminus \varphi_0\!\big(S(f)\cap U\big).
\]
Combining this identification with Theorem~\ref{thm:MSteich}, we obtain
\begin{equation} \label{eq:tech_identif}
\T(\hat U,f)\;\simeq\; \T\!\big((\GO(U)\cap \hat F(f))/f\big)\;\simeq\; \T\!\big(\varphi_0(U)\setminus \varphi_0(S(f)\cap U)\big),
\end{equation} 
and the set \(E:=(\GO(S(f))\cap U)/f\) corresponds, under this identification, to the set $$\tilde{E}:=\varphi_0(S(f)\cap U).$$
Thus, to finish the proof of Theorem~\ref{thm:teich}, it remains only to identify the uniformizing surface \(\varphi_0(U)\): in the contracting case, show \(\varphi_0(U)\simeq \C\) (so \(\T(\hat U,f)\simeq \T(\C\setminus \tilde{E})\)); in the non-contracting case, show \(\varphi_0(U)\simeq \D\) (so \(\T(\hat U,f)\simeq \T(\D\setminus \tilde{E})\)).

We start with the contracting case. 
Recall that $f$ is univalent in $W_n := \psi_n^{-1}(\D_r)$ for all $n$ and some fixed $r>0$. Since $\psi_n$ is a hyperbolic isometry,
$$W_n \;=\; D_{U_n}\!\big(f^n(z_0),\, R\big),
\quad R := 2\,\operatorname{arctanh}(r).$$ 
By Koebe's $1/4$-theorem, for all $n$, the set $\varphi_n(W_n)$ contains a disc of radius $r/4$ around the origin. For any $w\in U$, \eqref{eq:comp} and (the converse of) \eqref{eq:inf_prod} yield
\[
\varphi_m\!\big(f^m(w)\big)=\Big(\prod_{j=0}^{m-1}\lambda_j\Big)\varphi_0(w)\xrightarrow[m\to\infty]{}0,
\]
hence for some $m$ we have $\varphi_m(f^m(w))\in \D_{r/4}\subset \varphi_m(W_m)$; since the wandering domains are contracting, we have that $d_{U_m}(f^m(w),f^m(z_0))\to 0$ as $m\to\infty$, which implies that  $f^m(w)\in W_m$ for $m$ large enough. Thus
\[
U=\bigcup_{m\ge1} f^{-m}(W_m).
\]
Applying \eqref{eq:comp} once more $
\varphi_0 \big(f^{-m}(W_m)\big)\ \supset\ \frac{1}{\lambda_0\cdots \lambda_{m-1}}\D_{r/4},
$
so
\[
\varphi_0(U)\ \supset\ \bigcup_{m\ge1}\ \frac{1}{\lambda_0\cdots \lambda_{m-1}}\,\D_{\frac r4}=\C, 
\]
where the equality holds since  $\prod_{j=0}^{m-1}|\lambda_j|\to0$ by \eqref{eq:inf_prod}. Hence $\varphi_0(U)\simeq \C$. 

Moreover, in the contracting case, $\tilde E$ must be infinite.
If $\tilde E$ were finite, then 
\[
(\GO(U)\cap \hat F(f))/f \ \cong\ \varphi_0(U)\setminus \tilde E \ \cong\ \C\setminus\tilde E
\]
would be a surface of finite type; hence its Teichmüller space would be finite-dimensional. 
This contradicts Theorem~\ref{thm:main}, so $\tilde E$ cannot be finite, which makes it at least countably infinite.

Assume now that $U$ is semi-contracting. Recall from the proof of Theorem \ref{thm:simlin}, more precisely, \eqref{eq:lim5} and \eqref{eq:E2},  that
\[
\phi_0 \;=\; \lim_{m\to\infty}\frac{g_{m}\circ g_{m-1}\circ\cdots\circ g_{1}}{\lambda_{m}\lambda_{m-1}\cdots \lambda_{1}},
\]
which implies that the sequence $g_{m}\circ g_{m-1}\circ\cdots\circ g_{1}$ converges locally uniformly to some limit function $G:\D\to \D$ (since the denominator is not tending to 0). Under these conditions, and since $U$ is semicontracting, a theorem of Ferreira \cite[Theorem~1.1]{Fer23}  implies that $G$ is an inner function,  and therefore $\phi_0$ is a rescaling of an inner function.  Hence
$\varphi_0(U)\simeq\ \D \setminus O_0,$ for some (possibly empty) set $O_0$ of omitted values. Since any omitted value of an inner function is asymptotic, and asymptotic values are singular, we have that $O_0\subset \tilde{E}$, and so, by \eqref{eq:tech_identif} we are done.

Finally, assume that $U$ is eventually isometric: there exists $N$ such that
$f:U_k\to U_{k+1}$ is a hyperbolic isometry (hence a biholomorphism) for all $k\ge N$.
Fix $m\ge N$. Applying Lemma~\ref{lem:GOphi} at index $m$ yields a biholomorphism
\[
(\GO(U)\cap \hat F(f))/f \;\xrightarrow{\ \sim\ }\; \varphi_m(U_m)\setminus \varphi_m\!\big(S(f)\cap U_m\big).
\]
Since $f$ has no singular values on the isometric tail, $S(f)\cap U_m=\emptyset$, and therefore
\[
(\GO(U)\cap \hat F(f))/f \;\cong\; \varphi_m(U_m).
\]
But $U_m\simeq\D$, and $\varphi_m$ is a biholomorphism onto its image, so $\varphi_m(U_m)\simeq\D$.
Hence $\T(\hat U,f)\simeq \T(\D\setminus \tilde E)$, which concludes the proof. 

\section{Indiscrete characterisation: Proof of Theorem \ref{thm:indiscrete_intro}}\label{sec:proofindiscrete}
In this section we prove Theorem \ref{thm:indiscrete_intro}. Informally, we show that the orbit relation is indiscrete on a component of $\widehat F(f)$ in a wandering domain precisely when, after uniformization, the dynamics of $f$ on the forward orbit of that component is conjugate to iterated power maps on annuli (or punctured discs). This extends the picture already known for multiply connected wandering domains, see \cite[Corollary 6.9]{EFGP24} and \cite[Theorem 5.2]{BRS13}, and shows that the same geometric and dynamical rigidity occurs for components of $\widehat{F}$ lying inside simply connected wandering domains.

The proof separates into two parts. The geometric part establishes that \textit{indiscreteness} forces double connectivity of the components, Theorem \ref{thm:ind}. The dynamical part then follows from a Böttcher–type conjugacy for non-autonomous coverings, which yields the power-map dynamics; this is provided by results in \cite{Fer22} (see Lemma \ref{gust}). Recall that the \emph{connectivity} of a domain $U\subset\C$, denoted $c(U)$, is the number of connected components of $\widehat\C\setminus U$; we say that $U$ is \emph{doubly connected} if this number equals $2$.

\subsection{Tools from hyperbolic geometry: collars}
Our arguments have a strong hyperbolic flavour. For convenience, we gather here the collar results we use (collar theorem, thick–thin decomposition, and distance–injectivity estimates), taken from \cite[Chapters~2–3]{Hub06} and \cite{Thu97}. We refer the reader to these sources for details and definitions.

By an {\em isometry} between Riemann surfaces we mean a distance-preserving biholomorphic map.
By a {\em horocycle around a puncture $p$} in a hyperbolic Riemann surface $S$ we mean the image under a universal covering of a horocycle in $\D$ that is homotopic to the $\{p\}$.

\begin{lemma}[The collar theorem]\label{lem:collar}
	Let $S$ be a hyperbolic Riemann surface.
	\begin{itemize}
		\item If $\gamma\subset S$ is a simple closed geodesic of (hyperbolic) length $\ell$, then 
		\[ A_\gamma := \{z\in S\colon d_S(z, \gamma) < \eta(\ell)\}, \quad
		\text{ where }
		\quad \eta(\ell) = \frac{1}{2}\log\frac{\cosh(\ell/2) + 1}{\cosh(\ell/2) - 1}, \]
		is isometric to the $\eta(\ell)$-neighbourhood of the unique closed geodesic in an annulus of modulus $\pi/\ell$. In particular, $A_\gamma$ is a topological annulus,  known as the \emph{collar around~$\gamma$}.
		\item If $p$ is a puncture (cuspidal end) 
		of $S$, let $\gamma$ be a horocycle 
		of length $2$ around $p$, and denote by $A_p$ the subregion of $S$ whose boundary is $\gamma\cup\{p\}$. Then, $A_p$ is isometric to the region bounded by the horocycle of length $2$ in $\DD$, and it is known as the \emph{collar around $p$}.
	\end{itemize}
	\noindent Furthermore, collars around disjoint closed geodesics are disjoint, and collars around punctures are disjoint from all other collars on $S$.
\end{lemma}
Collars give us ``building blocks'' of Riemann surfaces whose geometry is very well-understood. 
\begin{defi}[Injectivity radius]
	Let $S$ be a hyperbolic Riemann surface, and let $p\in S$. The \textit{injectivity radius} of $S$ at $p$ is defined as
	\[ 
	\inj_S(p) := \frac{1}{2}\inf\left\{\ell_S(\gamma)\colon\text{$\gamma\subset S$ is a non-contractible loop through $p$}\right\}. 
	\]
\end{defi}

Equivalently, one can see that $\inj_S(p)$ is the largest radius of a hyperbolic disc around $p$ which is still a topological disc. 
Points with $\inj_S(p)\ge \varepsilon$ (the \emph{$\varepsilon$-thick part}) have a uniform 
disc-like neighbourhood of radius $\varepsilon$. When $\inj_S(p)$ is smaller than a universal 
threshold $M$, the point must lie in the \emph{thin part}, which is precisely a 
collar: either a cusp neighbourhood or the collar around a short simple closed geodesic. 
The next lemma makes this precise.

\begin{lemma}[Thick-thin decomposition]\label{lem:thickthin}
	There exists a universal constant $M > 0$ such that, for any hyperbolic Riemann surface $S$, any point $p\in S$ such that $\inj_S(p) < M$ belongs to either a collar around a puncture, or a collar around a simple closed geodesic.
\end{lemma}
Finally, we must delve further into the geometry of collars. 

\begin{lemma}\label{lem:smallcollars}
	Let $S$ be a hyperbolic Riemann surface, let $C\subset S$ be a collar (either a cusp collar around a puncture or a geodesic collar around a simple closed geodesic), and let $p\in C$.  
	Denote $d=d_S(p,\partial C)$. Then $\inj_S(p)\ \ge\ \frac{1}{2}\,e^{-d}.$
\end{lemma}

\begin{proof}
	We consider the two types of collars separately.
	
	Suppose first that $C$ is a cusp. We work in horocyclic coordinates $(p,t)$ on the standard cusp, where the metric is $dp^2+e^{2p}dt^2$ and a horocycle has length $e^{p}$ (see e.g. \cite[p.111, eq. (4.4.2)]{Bus}).
	Since the boundary horocycle of $C$ is chosen to have length $2$ , then it sits at $p_0=\log 2$. Then, $d=p_0-p$ in this case, and\ $\inj_S(p)=\tfrac12 e^{p}=\tfrac12 e^{p_0-d}=e^{-d}\ge \tfrac12 e^{-d}$.

	Suppose now that $C=\mathcal C(\beta)$ is the collar about a simple closed geodesic $\beta$, and denote $L=\cosh\!\big(\tfrac12\ell(\beta)\big)$. By \cite[Theorem 4.1.6(iii)]{Bus},
\[
\sinh \left(\inj_S(p)\right)=L\cosh d-\sinh d\geq \cosh d-\sinh d=e^{-d}.
\]
By monotonicity of $\operatorname{arsinh}$, and using that $\operatorname{arsinh}x\ge x/2$ for $0\le x\le 1$,
	\[
	\inj_S(p)\ =\ \operatorname{arsinh}\!\big(\sinh(\inj_S(p))\big)\ \ge\ \operatorname{arsinh}(e^{-d}) \ge \tfrac12 e^{-d},
	\]
	which concludes the proof of the lemma.
\end{proof}

\subsection{Proof of Theorem \ref{thm:indiscrete_intro}}
For the rest of the subsection, let $f\colon\Cx\to\Cx$ be an entire function with a wandering domain $U$, let $V$ be a component of $\widehat F(f)\cap U$ and denote $V_n=f^n(V)$, $n\geq 0$, so that, in particular, $V_0=V$.
\begin{thm}[Absorbing collars]\label{thm:collars}
Suppose that the orbit relation in $V$ is indiscrete. Then there exist $N\in\mathbb{N}$ and a sequence of collars $(C_n)_{n\geq N}$, $C_n\subset V_n$, such that
\begin{enumerate}[(a)]
		\item For every compact subset $K\subset V$, $f^n(K)\subset C_n$ for all sufficiently large $n$;
		\item $f(C_n)\subset C_{n+1}$.
	\end{enumerate}
\end{thm}

\begin{proof}
Let $z\in V$. Since the orbit relation in $V$ is indiscrete, it follows from Theorem~\ref{thm:GOrelations}, see also \cite[Prop. 6.7]{EFGP24} for details, that $z$ is an accumulation point  of $\GO(z)$. We claim that then, $\inj_{V_n}(f^n(z))\to 0$ as $n\to\infty$. Indeed, choose $(z_n)\subset \GO(z)\cap V$ such that $z_n\to z$ and join each $z_n$ to $z$ by a distance-minimising geodesic arc $\gamma_n$; since these are points in the grand orbit of $z$, we can (by reshuffling the labels if necessary) say that 
$f^n(z_n) = f^n(z)$ for all $n$. Since $z_n\to z$, we also know that $\ell_V(\gamma_n)\to 0$.
Now, by the unique path-lifting property of covering maps, we know that $f^n(\gamma_n)$ is a closed, non-contractible geodesic loop in $V_n$ based at $f^n(z_n)=f^n(z)$. Since $f^n$ is also a local hyperbolic isometry, we see that
	\[ \inj_{V_n}(f^n(z)) \leq \frac{1}{2}\ell_{V_n}(f^n(\gamma_n)) = \frac{1}{2}\ell_V(\gamma_n)\to 0 \quad \text{ as $n\to\infty$,} \]
as claimed.


 Next, we show that for any compact subset $K\subset V$ there exist collars $C_n\subset V_n$ with $f^n(K)\subset C_n$ for all sufficiently large $n$. Fix $z\in V$. Since $\inj_{V_n}(f^n(z))\to 0$ as $n\to\infty$, Lemma~\ref{lem:thickthin} yields $N\in\N$ such that, for each $n\ge N$, there is a collar $C_n\subset V_n$ containing $f^n(z)$. Set $D_n:=d_{V_n}(f^n(z),\partial C_n)$, and let $K\subset V$ be a compact set. Choose $L>0$ with $d_V(w,z)\le L$ for all $w\in K$.  By Lemma~\ref{lem:smallcollars} we have $D_n\to\infty$. By Schwarz--Pick and the triangle inequality,
\[
d_{V_n}\!\big(f^n(w),f^n(z)\big)\le d_V(w,z)\le L\quad\text{for all }w\in K.
\]
Fix $n$ so large that $D_n>L$. Suppose, towards a contradiction, that there exists $w\in K$ with $f^n(w)\notin C_n$. Let $\beta_n$ be a distance-minimising geodesic in $V_n$ from $f^n(z)$ to $f^n(w)$. Since $f^n(z)\in C_n$ and $f^n(w)\notin C_n$, the curve $\beta_n$ meets $\partial C_n$ at some point $b_n$, and hence
\[
D_n=d_{V_n}(f^n(z),\partial C_n)\le d_{V_n}\!\big(f^n(z),b_n\big)\le \length(\beta_n)=d_{V_n}\!\big(f^n(z),f^n(w)\big)\le L,
\]
a contradiction to $D_n>L$. Therefore $f^n(w)\in C_n$ for all $w\in K$ once $n$ is large; that is, $f^n(K)\subset C_n$ for all large $n$.

Fix $z\in V$ and let $(C_n)$ be the collars constructed from its orbit. We claim that the sequence of collars is unique: for any other sequence $(C_n')$ constructed from $z'\in V$, every compact $K\Subset V$ eventually intersects both $C_n$ and $C_n'$, and since distinct collars are pairwise disjoint (Theorem~\ref{thm:collars}), it follows that $C_n=C_n'$ for all sufficiently large $n$.

Finally, we are left to show that $f(C_n)\subseteq C_{n+1}$. To see this, consider a point $w\in C_n$ and a length-minimising closed geodesic loop $\gamma$ through $w$. Since $f\colon V_n\to V_{n+1}$ is a covering, it induces an injective morphism from $\pi_1(V_n)$ into $\pi_1(V_{n+1})$, meaning that $f(\gamma)$ is a non-trivial loop through $f(w)$; by the Schwarz--Pick lemma, $\ell_{V_{n+1}}(f(\gamma)) \leq \ell_{V_n}(\gamma)$, and so $\inj_{V_{n+1}}(f(w)) \leq \inj_{V_n}(w)$. By Lemma \ref{lem:thickthin}, $f(w)$ belongs to $C_{n+1}$, which concludes the proof of the theorem.
\end{proof}

We have just seen in the previous theorem that when the orbit relation is indiscrete, there exist \textit{absorbing} collars that propagate forward under the dynamics and push the iterates into them. This leaves only two possible geometric outcomes, stated in the next theorem.

\begin{thm}[Geometry for indiscrete GO relations]\label{thm:ind}The grand orbit relation in $V$ is indiscrete if and only if one of the following holds.
	\begin{enumerate}[(i)]
		\item $V$ is an annulus of finite modulus and $f^n\colon V\to V_n$ is a covering map between annuli of degree $D_n$, with $D_n\to\infty$ as $n\to\infty$; or
		\item $V$, and each $V_n$, is a punctured disc, and the puncture is a critical point of $f$ for infinitely many values of $n$.
	\end{enumerate}
\end{thm}
\begin{proof}
It is easy to see that either (i) or (ii) imply that the grand orbit relation in $V$ is indiscrete. Indeed, let $\gamma$ be any noncontractible simple closed curve in $V$, and set $z\in\gamma$. Since $f^n|_V$ is a covering of increasing degree, the set  $f^{-n}(f^n(z)) \cap \gamma \subset \GO(z)\cap V$ has infinitely many points, which must therefore have an accumulation point on $\gamma\subset V$. Hence $\GO(z)$ is not a discrete set which implies that the grand orbit relation in $V$ is indiscrete, by Theorem \ref{thm:GOrelations}.
	
We now show the converse; assume that the orbit relation in $V$ is indiscrete. Let us start with the following observation:
\begin{lemma}
If $V$ is doubly connected, then the same is true of $V_n$ for all $n\geq 1$.
\end{lemma}
\begin{proof}
Consider the covering map $f^n\colon V\to V_{n+1}$. If it has finite degree, the Riemann--Hurwitz formula (see e.g. the form given in \cite[Theorem 3(i)]{Bol99}) tells us that
	$$ c(V) - 2 = D_n\cdot\left(c(V_n) - 2\right) + \delta_V(f^n), $$
where $c(V)$ denotes the connectivity of the domain $V$,  $D_n = \deg(f^n|V)$ and $\delta_V(f^n)$ is the number of critical points of $f^n$ in $V$, counted with multiplicity. Since it is a covering, the latter is zero, and so therefore $c(V) = 2$ implies $c(V_n) = 2$. If, on the other hand, $f^n|_V$ has infinite degree, we can (since $f$ is entire, and thus meromorphic) invoke \cite[Theorem~3(ii)]{Bol99}: if $c(V_n) \geq 3$, we have $c(V) = \infty$, which is a contradiction. However, we also cannot have $c(V_n) = 1$, for $\pi_1(V)\simeq \mathbb{Z}$ injects into $\pi_1(V_n)$ via $f^n_*$. This completes the proof.
\end{proof}

With that done, notice that the maps also cannot be eventually degree-one isometries on $(V_n)$; hence
	$d_n=\deg\!\big(f|_{V_n}:V_n\to V_{n+1}\big)\ge 2$ for infinitely many $n$, and thus
	$D_n=d_{n-1}\cdots d_0\to\infty$ (in the punctured-disc case: the local degree at the puncture is $\ge 2$ for infinitely many $n$).
	
	Now, let $(C_n)_{n\geq N}$ be the sequence of absorbing collars given by Theorem \ref{thm:collars}. For the same reason, since the restriction $f\vert_{C_n}\colon C_n \to C_{n+1}$ is a covering it preserves the type of component: either $C_n$ is a  collar around a geodesic $\sigma_n$ for all $n\geq N$, or is a  collar around a puncture for all $n\geq N$. (Indeed, a covering sends nontrivial closed geodesics to nontrivial closed geodesics, so a geodesic collar cannot map into a cusp; conversely, a cusp contains no closed geodesics, so a cusp collar cannot map onto a geodesic collar.) We divide the proof in two cases, based on this dichotomy.
	
	Assume first that the sets $C_n$ are collars around punctures. If the connectivity of  $V$ is more than $2$, then $V$ contains a closed geodesic $\gamma$; since $\gamma$ is a compact subset of $V$, by Theorem~\ref{thm:collars} there exists $n\in\mathbb{N}$ such that $f^n(\gamma)\subset C_n$ for all $n$ large enough. Since $f^n\colon V\to V_n$ is a covering, $f^n(\gamma)$ is a closed geodesic contained in the collar around a puncture -- a contradiction. Hence $V$ is at most doubly connected. Moreover, for large $n$ the sets $C_n\subset V_n$ are cusp collars, so $V_n$ has a puncture for $n$ large enough; since coverings lift punctures to punctures, we deduce that $V_j$ is a punctured disc for all $j\geq 0$.

	Assume now that, for all large $n$, the sets $C_n$ are collars around simple closed geodesics $\sigma_n$. Fix $N$ sufficiently large and let $\sigma\subset V$ be a component of $f^{-N}(\sigma_N)\cap V$; then $\sigma$ is a simple closed geodesic in $V$. Suppose, for a contradiction, that there exists another simple closed geodesic $\gamma\subset V$ with $f^N(\gamma)=\sigma_N$ and such that $\gamma\cap\sigma=\varnothing$. Let $n\geq N$, $z\in\sigma$ and $w\in\gamma$ with $f^n(z)=f^n(w)\in\sigma_n$. Let $\alpha$ be a distance-minimising geodesic segment in $V$ joining $z$ to $w$; by uniqueness of the geodesic flow, $\alpha$ is transversal to both $\sigma$ and $\gamma$. Since $f^n(z)=f^n(w)$, the image $f^n(\alpha)$ is a closed geodesic \emph{arc}\footnote{There is a subtle distinction between a \textit{closed geodesic} (which is smooth) and a \textit{closed geodesic arc} (which is continuous, but not smooth).} based at $f^n(z)$. For $n$ large, $f^n(\alpha)\subset C_n$ by Theorem~\ref{thm:collars}. Because $f$ is angle-preserving, $f^n(\alpha)$ is transverse to $\sigma_n$ at $f^n(z)$, hence $f^n(\alpha)\neq\sigma_n$. This contradicts Lemma~\ref{lem:collar}, which asserts that in a geodesic collar the only closed geodesic arc based at a point of $\sigma_n$ is $\sigma_n$ itself. Therefore $V$ contains at most one simple closed geodesic. It follows that $V$ is an annulus, and hence (since $f^n|_V$ is a covering) each $V_n$ is also an annulus.
\end{proof}
Finally, we shall use the following version of \cite[Lemma~2.4]{Fer22}.
\begin{lemma}\label{gust} Suppose the domains $V_n$ are doubly-connected for all $n\geq 0$. Then there exists a sequence of conformal maps $\phi_n\colon V_n\to A_n$, where $A_n$ is either an annulus of finite modulus or a punctured disc, such that 
	\[\phi_{n+1}\circ f(z) =  (\phi_n(z))^{d_n} \quad \text {for all } z\in V_n.\]
\end{lemma}	
\begin{remark} The original statement in \cite[Lemma~2.4]{Fer22}
concerns (full) wandering domains rather than components of $\widehat{F}(f)$ in wandering domains. However, its proof only depends on the map $f$ being a covering from the component to its image, and so it also applies to the sets $V_n$. Moreover, in the original statement, $A_n$ is necessarily an annulus, since (full) wandering domains cannot be punctured discs given that the Julia set is a perfect set. However, the proof can be adapted to the punctured disc case. 
\end{remark} 

Theorem~\ref{thm:indiscrete_intro} follows immediately from Theorem~\ref{thm:ind} together with Lemma~\ref{gust}.

\bibliographystyle{alpha}

\bibliography{ref}

\end{document}